\documentclass[12pt]{article}
\usepackage{amsmath,amssymb}

\newcommand{\ncm}{\newcommand}
 \ncm{\R}{\mathbb{R}}
 \ncm{\C}{\mathbb{C}}
 \ncm{\Q}{\mathbb{Q}}
 \ncm{\N}{\mathbb{N}}
 \ncm{\Ad}{\mbox{\rm Ad}}
 \ncm{\ad}{{\rm ad}}
 \ncm{\Aut}{\mbox{\rm Aut}}
 \ncm{\Sp}{{\rm Sp}}
 \ncm{\Spec}{{\rm Spec}}
 \ncm{\supp}{{\rm supp}}
 \ncm{\Tr}{{\rm Tr}}
 \ncm{\id}{{\rm id}}
 \ncm{\ra}{\rightarrow}
 \ncm{\cstar}{C$^*$-algebra}
 \ncm{\M}{{\mathcal{M}}}
 \ncm{\Ex}{{\mathcal Ex}}
 \ncm{\E}{{\mathcal E}}
 \ncm{\G}{{\mathcal G}}
 \ncm{\F}{{\mathcal F}}
 \ncm{\D}{{\mathcal D}}
 \ncm{\T}{\mathbb{T}}
 \ncm{\V}{{\mathcal V}}
 \ncm{\K}{{\mathcal K}}
 \ncm{\CO}{{\mathcal O}}
 \ncm{\I}{{\mathcal I}}
 \ncm{\U}{{\mathcal U}}
 \ncm{\Hil}{{\mathcal H}}

 \ncm{\Z}{{\mathbb{Z}}}
 \ncm{\eps}{\epsilon}
 \ncm{\ran}{{\rangle}}
 \ncm{\lan}{{\langle}}

 \newtheorem{theo}{Theorem}[section]
\newtheorem{cor}[theo]{Corollary}
\newtheorem{lem}[theo]{Lemma}
\newtheorem{prop}[theo]{Proposition}
\newtheorem{remark}[theo]{Remark}
\newtheorem{definition}[theo]{Definition}
\newtheorem{example}[theo]{Example}

\newenvironment{pf}{{\it Proof.}}{QED\vspace{3mm}}

\title{Approximately inner flows}
\author{A. Kishimoto\footnote{ E-mail: akiksmt@r3.ucom.ne.jp; Retired from Hokkaido University}}
\begin{document}

\maketitle

\begin{abstract}
When $\alpha$ is an approximately inner flow on a C$^*$-algebra $A$ and commutes with an automorphism $\gamma$ of $A$ we may extend $\alpha$ to a flow $\bar{\alpha}$ on the crossed product $A\times_\gamma\Z$ by setting $\bar{\alpha}_t(U)=U$ where $U$ is the canonical unitary implementing $\gamma$ in $A\times_\gamma\Z$ and ask whether $\bar{\alpha}$ is also approximately inner or not. We will consider very specific examples of this type; some of which we can answer affirmatively.
\end{abstract}

\section{Introduction}

Let $(k_n)$ be a sequence in $\N$ and $\beta_n$ a flow on the C$^*$-algebra $M_{k_n}$ of $k_n\times k_n$ matrices. We define a one-parameter group $\beta$ of automorphisms of $B=\prod_n M_{k_n}$ by $\beta_t=\prod_n \beta_{n,t},\ t\in\R$. Let $B_\beta$ denote the maximal C$^*$-subalgebra of $B$ on which $t\mapsto \beta_t(x)$ is continuous; so $\beta$ is a strongly continuous one-parameter group of automorphisms of $B_\beta$ or a {\em flow} on $B_\beta$. Since $I=\bigoplus_n M_{k_n}$ consists of $x\in B$ converging to zero, it follows that $I$ is a $\beta$-invariant ideal of $B$ with $I \subset B_\beta$. We also denote by $\beta$ the flow on the quotient $B_\beta/I$ induced by $\beta$. We recall that  a flow $\alpha$ on a separable C$^*$-algebra $A$ is called an MF flow if there is an embedding $\phi$ of $A$ into $B_\beta/I$ such that $\phi\alpha_t=\beta_t\phi$ for some sequence $(k_n)$ in $\N$ and $(\beta_{n})$ \cite{QD2}. Suppose that $A$ is an MF C$^*$-algebra (i.e., embeddable into $B/I$ for some sequence $(k_n)$ in $\N$). Then if $\alpha$ is approximately inner it follows that $\alpha$ is an MF flow (see \cite{AI-MF}). A question remains whether the converse holds for some class of MF C$^*$-algebras $A$, e.g., AF C$^*$-algebras. (Without any condition on $A$ this does not hold, e.g., there is a non-trivial MF flow on a commutative C$^*$-algebra.)

Since $B_\beta/I$ is inseparable, $\beta$ is not an MF flow on $B_\beta/I$ by definition. We will give a necessary and sufficient condition for $\beta$ to be  approximately inner. Inspired by its proof we give an example of an MF flow which is not AI; unfortunately the C$^*$-algebra is not simple, rather similar to the Toeplitz algebra. See Section 2.

Then we turn to specific examples of MF flows which may not be approximately inner (AI for short). However we could not prove that they include non-AI flows; instead we showed that some of them are actually AI, whose proofs we hope are somewhat non-trivial; that is why we are presenting them here. We shall now describe those examples in detail.

Let $A$ be a unital C$^*$-algebra $A$. If $\alpha$ is a flow on $A$ and $\gamma$ is an automorphism of $A$ such that $\alpha_t\gamma=\gamma\alpha_t,\ t\in\R$, we define a flow $\bar{\alpha}$ on the crossed product $A\times_\gamma\Z$ such that $\bar{\alpha}_t|A=\alpha_t$ and $\bar{\alpha}_t(U)=U$, where $U$ is the canonical unitary of $A\times_\gamma\Z$ implementing $\gamma$ on $A$.

We are concerned with the following problem: If $\alpha$ is AI and $\gamma$ is not so wild, can we conclude that $\bar{\alpha}$ is approximately inner? (This problem was first taken up in \cite{AI-MF}.) Without any condition on $(A,\gamma)$ it does not follow that $\bar{\alpha}$ is AI. To show this one may use the fact that the Cuntz algebra $\CO_2$ can be expressed as a crossed product of $M_{2^\infty}$ by a trace-scaling endomorphism \cite{C77}. More precisely let $A=M_{2^\infty}\otimes\K$ with $\K$ the C$^*$-algebra of compact operators and let $\gamma$ be an automorphism of $A$ such that $\tau\gamma=2\tau$ where $\tau$ is a non-trivial lower semi-continuous trace; then $\gamma$ is unique up to inner automorphisms \cite{EEK98} and $A\times_\gamma\Z$ is isomorphic to $\CO_2\otimes\K$.

More specifically we denote by $\Hil$ the Hilbert space spanned by a complete orthonormal family $(\xi_\Lambda)$ of vectors indexed by all finite subsets $\Lambda$ of $\N$ and define a unitary flow $U$ by $U_t\xi_\Lambda=e^{ip|\Lambda|t}$ for some $p\in\R$ where $|\Lambda|$ is the cardinality of $\Lambda$. We define a flow on $A=M_{2^\infty}\otimes \K(\Hil)$ by
$$
\alpha_t=\bigotimes_n \Ad(1\oplus e^{ipt})\otimes \Ad\,U_t.
$$
Let $S$ denote the shift on $\N$, i.e., $S(n)=n+1$, and $V$ denote the isometry on $\Hil$ defined by $V\xi_\Lambda=\xi_{S(\Lambda)}$, where $S(\Lambda)=\{S(\lambda)\ |\ \lambda\in \Lambda\}$. Let $e_{i,j},\ i,j=1,2$ be matrix units on $\Hil$ such that $e_{1,2}+e_{2,2}=1$ and $e_{2,2}\xi_\Lambda=\xi_\Lambda$ if $1\in \Lambda$; $=0$ otherwise. We define an automorphism $\gamma$ of $A=M_{2^\infty}\otimes\K(\Hil)$ by requiring that $\gamma$ on $M_{2^\infty}=\bigotimes_{n=0}^{-\infty} M_2$ is the shift to the right sending the last factor $M_2$ at $0$ onto the algebra $M_2$ generated by $e_{i,j}$ on $\Hil$ by $e_{i,j}\otimes 1_\Hil\mapsto 1\otimes e_{i,j}$ and $\gamma$ on $\K(\Hil)$ is the endomorphism given by $x\mapsto \sum_{i=1}^2 e_{i,1}VxV^*e_{1,i}$. An inspection shows that $\alpha_t$ and $\gamma$ commutes with each other. Then the flow $\bar{\alpha}$ on $A \times_\gamma\Z$ identifies with an extension of the quasi-free flow on $\CO_2$ defined by $s_1\mapsto s_1$ and $s_2\mapsto e^{ipt}s_2$, which is not AI because it does not have a KMS state. This implies that $\bar{\alpha}$ is not AI. (The C$^*$-algebra $\CO_2\otimes\K$
is purely infinite, far from MF C$^*$-algebras; so it is still desirable to have such an example with a more benign $\gamma$.)

Avoiding such a situation we consider the following triple $(A,\alpha,\gamma)$. Let $F$ be a finite-dimensional C$^*$-algebra and let $A(m)=F$ for each $m\in\Z$ and define $A=\bigotimes_{m\in\Z}A(m)$. For a finite subset $\Lambda$ of $\Z$ let $A(\Lambda)$ denote the C$^*$-subalgebra of $A$ given by $\bigotimes_{m\in \Lambda}A(m)$. We identify $A(\{m\})$ with $A(m)$.

Let $\gamma$ denote the translation automorphism of $A$; $\gamma(A(m))=A(m+1)$. We will also regard $\gamma$ as an action of $\Z$ on $A$ by setting $\gamma_n=\gamma^n,\ n\in\Z$. Let $\Phi$ be a $\gamma$-invariant potential; i.e., $\Phi$ is a function of the non-empty finite subsets of $\Z$ into $A_{sa}$ such that $\Phi(\Lambda)\in A(\Lambda)$, $\gamma(\Phi(\Lambda))=\Phi(\Lambda+1)$ for each finite subset $\Lambda\not=\emptyset$, and $\sum_{ \Lambda \ni 0} \|\Phi(\Lambda)\|<\infty$. For a finite subset $\Delta$ of $\Z$ we set
$$
H(\Delta)=\sum_{\Lambda\subset \Delta}\Phi(\Lambda).
$$
If $\Phi$ satisfies
$$
\|\Phi\|_\lambda\equiv\sum_{\Lambda\ni0} e^{\lambda|\Lambda|}\|\Phi(\Lambda)\|<\infty
$$
for some $\lambda>0$ (or $\Phi(\Lambda),\ \Lambda\subset \Z^d$ mutually commute) then $\Ad\,e^{itH(\Lambda)}$ converges and defines a flow $\alpha^\Phi$ on $A$ as $\Lambda$ increases to $\Z$ \cite{BR2}. We then define a flow $\bar{\alpha}^\Phi$ on $A\times_\gamma\Z$ by setting $\bar{\alpha}^\Phi_t|A=\alpha^\Phi_t$ and $\bar{\alpha}^\Phi_t(U)=U$, where $U$ is the canonical unitary implementing $\gamma$.

We do not know whether $\bar{\alpha}^\Phi$ is AI or not for a general $F$. But we will show that $\bar{\alpha}^\Phi$ is AI if $F$ is a full matrix algebra or $F$ is of the form $M_k\otimes F_1$ with $k\geq 2$ where $F_1$ a finite-dimensional C$^*$-algebra. This answers a problem left out in \cite{AI-MF}, where $\bar{\alpha}^\Phi$ was shown to be quasi-diagonal (or QD) at least when $F$ is a full matrix algebra but was not decided to be AI or not. See \cite{QD1,QD2,QD3} for more on QD flows. See Section 3 for the full matrix algebra case and Section 4 for the other case.

We recall another condition on flows. The flow $\alpha$ on $A$ is called {\em continuously AI} if there is a continuous function $h:[0,\infty)\ra A_{sa}$ such that $\alpha_t(x)=\lim_{s\ra\infty}\Ad\,e^{ith(s)}(x), x\in A$ uniformly in $t\in [-1,1]$ \cite{K06}. The $\alpha^{\Phi}$ defined above based on a potential $\Phi$ is not only AI but also continuously AI. (The generator $\delta^{\Phi}$ of $\alpha^{\Phi}$ has, as a core, $\bigcup_{\Lambda}A(\Lambda)$, where $\Lambda$ runs over all finite subsets of $\Z$. In this case if $(h_n)$ in $A_{sa}$ satisfies $\alpha_t=\lim_n \Ad\,e^{ith_n}$ then the linear extension of $h:\{0,1,2,\ldots\}\ra A_{sa}$ defined by $h(n)=h_n$ to a function on $[0,\infty)$ will automatically satisfy $\alpha_t=\lim_s\Ad\,e^{ith(s)}$; see \cite{Sak}.) I do not know whether being continuously AI is strictly stronger than being AI and was hoping to check whether $\bar{\alpha}^\Phi$ is continuously AI or not. We will show that if there is no interaction (i.e., $\Phi(\Lambda)=0$ for all $\Lambda$ except for singletons) then $\bar{\alpha}^\Phi$ is continuously AI for a matrix algebra $F$ but leave the problem undecided of whether there is $\Phi$ and $F$ such that $\bar{\alpha}^\Phi$ is not continuously AI. See Section 5.

The condition of being continuously AI was introduced in connection of the lifting problem: If $\beta$ is an AI flow on $B=A/I$ with $I$ an ideal of $A$, then is there a flow $\alpha$ on $A$ which induces $\beta$ on the quotient $A/I$? This is affirmative if $\beta$ is continuously AI. (If we put another condition on $\alpha$ that $\alpha|I$ should be universally weakly inner, this is necessary, for this statement to hold for any extension of $B$, at least when $B$ is simple. Consider the extension $A=C[0,1]\otimes B$ with $I=C_0[0,1)\otimes B$.)  We tried to elucidate the situation better to no avail.

We will conclude this note by giving some excuses again in Section 6, where it is also shown that the above candidates for AI flows are all MF (and QD) flows.

%%%%%%%%%%%%%%%%%%%%%%%%%%%%%%%%%%%%%%%%%%%%%%%%%%%%%%%%%%%%%%%%%%%
\section{The flow $\beta$ on $B_\beta/I$}

 We recall $B=\prod_n M_{k_n}$, $I=\bigoplus_n M_{k_n}$, and $\beta_t=\prod_n\beta_{n,t}$ for some $(k_n)$ and $(\beta_n)$ with $\beta_n$ a flow on $M_{k_n}$. The C$^*$-subalgebra $B_\beta$ of $B$ is defined as the maximal $\beta$-invariant C$^*$-subalgebra on which $t\mapsto \beta_t$ is continuous. We first recall the following result concerning KMS states \cite{QD2}.

\begin{prop}\label{KMS}
Let $B=\prod_{n=1}^\infty M_{k_n}$ and $I=\bigoplus_{n=1}^\infty M_{k_n}$ for some sequence $(k_n)$ in $\N$. Let $\beta_n$ be a flow on $M_{k_n}$ and let $\beta_t=\prod_n \beta_{n,t}$. Then the flow $\beta$ on $B_\beta/I$ has KMS states for all inverse temperatures.
\end{prop}
\begin{pf}
Fix an inverse temperature. Each $\beta_n$ has a unique KMS state $\omega_n$ on $M_{k_n}$. Extend $\omega_n$ to a state $\bar{\omega}_n$ of $B_\beta$ by $\bar{\omega}_n(x)=\omega_n(x_n)$ for $x=(x_n)\in B_\beta$, which is a KMS state with respect to $\beta$. Let $\omega$ be an accumulation point of $(\bar{\omega}_n)$, which is a KMS state such that $\omega(I)=0$. Hence one can regard $\omega$ as a state of $B_\beta/I$. Thus $\beta$ has a KMS state on $B_\beta/I$ for any inverse temperature.
\end{pf}

Let $h_n\in (M_{k_n})_{sa}$ be such that $\beta_{n,t}=\Ad\,e^{ith_n}$ and $\Spec(h_n)\cap (-\infty,0]=\{0\}$ and let $u_t=(e^{ith_n})_n\in B_\beta/I$. Then $\beta_t=\Ad\,u_t$, i.e., $\beta_t$ is an inner automorphism for each $t$. Unless $(h_n)$ is bounded, $t\mapsto u_t$ is not continuous.

\begin{prop}
In the above situation for each $c>0$ let
$$
\ell_n(c)=\max\{b-a-c\ |\ [a,b]\subset \Spec(h_n)+[0,c]\}.
$$
Then $\beta$ is uniformly continuous on $B_\beta/I$ (i.e., $\delta_\beta$ is bounded) if and only if
$$
\sup_{c>0} \limsup_n\ell_n(c)<\infty.
$$
\end{prop}
\begin{pf}
Suppose $\sup_{c>0}\limsup_n \ell_n(c)=\infty$. For any $M>0$ there is a $c>0$ such that $\limsup_n\ell_n(c)>M$. Then there is an infinite subset $J$ of $\N$ such that $\ell_n(c)>M$ for $n\in J$. There is a pair $a_n,b_n\in \Spec(h_n)$ for $n\in J$ such that $M<b_n-a_n\leq M+c$ and hence there is a partial isometry $u_n\in M_{k_n}$  such that $\beta_{n,t}(u_n)=e^{it(b_n-a_n)}u_n$. Set $u_n=0$ for $n\not\in J$ and let $u=(u_n)\in B_\beta/I$, which is in the domain of $\D(\delta_\beta)$ and satisfies $\|\delta_\beta(u)\|>M$. Therefore $\beta$ is not uniformly continuous.

Suppose $0<g=\sup_{c>0}\limsup_n \ell_n(c)<\infty$. If $\Spec(\beta)\cap (g,\infty)\not=\emptyset$ then there is $x\in B_\beta/I$ such that $\Spec_\beta(x)\subset [a,b]\subset (g, \infty)$ for some $a,b>0$. Then it follows that $\Spec(\beta_n)\cap [a,b]\not=\emptyset$ for a infinitely many $n$, which implies that $\ell_n(b)\geq a$ for such $n$ or $\limsup_n\ell_n(b)\geq a>g$. This contradiction shows that $\Spec(\beta)\subset [-g,g]$; thus $\beta$ is uniformly continuous.
\end{pf}

Note that $(\|h_n\|)$ need not be bounded for $\beta$ to be uniformly continuous.

\begin{prop} \label{Toeplitz}
In the above situation $\beta$ is approximately inner on $B_\beta/I$ if and only if
$(\ell_n(c))$ is bounded for all $c>0$.

\end{prop}
\begin{pf}
Suppose that $\ell_n(c)\ra\infty$ as $n\ra\infty$ for some $c>0$. We choose $[a_n,b_n]\subset \Spec(h_n)+[0,c]$ such that $\ell(c)+c=b_n-a_n$. Let $p_{n,1},p_{n,2},\ldots,p_{n,s_n}$ be an increasing sequence of eigenvalues of $h_n$  such that $p_{n,1}=a_n$, $c\leq p_{n,i}-p_{n,i-1}< 2c$, and $b_n-2c< p_{n,s_n}(\leq b_n-c)$. Let $e^{(n)}_{i,j},\ i,j=1,2,\ldots,s_n$, be a family of matrix units such that $h_ne^{(n)}_{i,i}=p_{n,i}e^{(n)}_{i,i}$. We set
$$
u_n=\sum_{i=2}^{s_n} e_{i,i-1}^{(n)}
$$
which is a partial isometry such that $\beta_{n,t}(u_n)=\sum_{i=2}^{s_n}e^{it(p_{n,i}-p_{n,i-1})}e^{(n)}_{i,i-1}$. Since $\Spec_{\beta_n}(u_n)\subset [c,2c]$, it follows that $u=(u_n)\in B_\beta$.

Suppose that $\beta$ is approximately inner. Then there is a net $(b_k)$ in $(B_\beta/I)_{sa}$ such that $\delta_\beta$ is the graph limit of $\ad\,ib_k$ on $B_\beta/I$. Since for any finite-dimensional subalgeba $F$ of $\D(\delta_\beta)$ one may suppose that $\ad\,ib_k|F=\delta_\beta|F$, we impose the condition that $\ad\,ib_k=0$ on the three-dimensional subalgebra $\{uu^*,u^*u\}$. Since we are concerned about only $u$ we assume that $(b_k)$ is a sequence. There is a sequence $z_k$ in $B_\beta/I$ such that $z_k\ra u$ and $\ad\,ib_k(z_k)\ra \delta_\beta(u)$. By replacing $z_k$ by $uu^*z_ku^*u$ and functional calculus we may suppose that $z_k$ is a partial isometry such that $z_kz_k^*=uu^*$ and $z_k^*z_k=u^*u$.  Let $(b_{k,n})_n$ (resp. $(z_{k,n})$) be a representative of $b_k$ (resp. $z_k$) such that $b_{k,n}$ is self-adjoint and satisfies $\ad\,ib_{k,n}|\{u_nu_n^*,u_n^*u_n\}=0$ (resp. $z_{k,n}$ a partial isometry such that $z_{k,n}z_{k,n}^*=u_nu_n^*$ and $z_{k,n}^*z_{k,n}=u_n^*u_n$). Thus we deduce
$$
\limsup_n\|[b_{k,n},z_{k,n}]-[h_n,u_n]\|=\|[b_k,z_k]+i\delta_\beta(u)\|.
$$
Since $\|z_{k,n}^*[b_{k,n},z_{k,n}]-u_n^*[h_n,u_n]\|\leq \|[b_{k,n},z_{k,n}]-[h_n,u_n]\|+\|z_{k,n}-u_n\|\|[h_n,u_n]\|$, one finds, for any $\eps>0$, $k$ and $n$ such that
$$
\|z_{k,n}^*b_{k,n}z_{k,n}-z_{k,n}^*z_{k,n}b_{k,n} -(u_n^*h_nu_n-u_n^*u_nh_n)\|<\eps.
$$
Taking the trace on $u_n^*u_n\vee u_nu_n^*$ we obtain that
$$
|\Tr(b_{k,n}u_nu_n^*-b_{k,n}u_n^*u_n)-\Tr(h_nu_nu_n^*-h_nu_n^*u_n)|\leq \eps \Tr(u_n^*u_n\vee u_nu_n^*),
$$
which implies that
$$
p_{n,s_n}-p_{n,1}-2\|b_{k,n}\|\leq \eps (\frac{p_{n,s_n}-p_{n,1}}{c}+1)
$$
as the rank of $u_n^*u_n\vee u_nu_n^*$ is at most $(p_{n,s_n}-p_{n,1})/c +1$.
Since $\limsup_n\|b_{k,n}\|=\|b_k\|$ and $\lim_n(p_{n,s_n}-p_{n,1})=\infty$, this is a contradiction for $\eps<c$. Hence $\beta$ is not AI. If $\limsup_n\ell_n(c)=\infty$ there is a subsequence $(k_n)$ in $\N$ such that $\lim_n \ell_{k_n}(c)=\infty$ and we reach the same conclusion.

On the contrary suppose that $(\ell_n(c))$ is bounded for all $c>0$. Let $h_n=\sum_i\lambda_ip_i$, where $(p_i)$ is a family of rank-one projections with $\sum_ip_i=1$. Then the family $\{\lambda_i\}$ is divided into a few groups such that each group has values in an interval whose length is at most $\ell_n(c)$ and the neighboring intervals are separated by more than $c$. We then define $h_n(c)$ in the form $\sum_i\lambda_i'p_i$  by translating the eigenvalues of $h_n$ in each group by the same value so that the norm of $h_n(c)$ is at most $\ell_n(c)$. Set $h(c)=(h_n(c))\in B_\beta$. Then it follows that $\delta_\beta(x)=\ad\,ih(c)(x)$ for $x\in B_\beta/I$ with $\Spec_\beta(x)\subset (-c,c)$. This implies that $\delta_\beta$ is the graph limit of $\ad\,ih(c)$ as $c\ra\infty$. This concludes the proof that $\beta$ is approximately inner.
\end{pf}

We recall the following implications for flows on separable unital MF C$^*$-algebras: AI $\Rightarrow$ MF $\Rightarrow$ KMS; where KMS means the existence of KMS states for all inverse temperatures, which follows from Proposition \ref{KMS}. Whether the converse of each implication holds for simple C$^*$-algebras is not known.

We give two examples, which are taken from the essence of the proof of Proposition \ref{Toeplitz}. The second one is more pertinent (though the C$^*$-algebra is not even prime); the first one is a simpler version.

\begin{example}
Let $\mathcal{T}=C^*(S)$ be a Toeplitz algebra, where $S$ is the one-sided shift on a Hilbert space: $S\xi_i=\xi_{i+1}$ where $\xi_0,\xi_1,\xi_2,\ldots$ is a complete orthnormal system of $\Hil$. Define a flow $\alpha$ on $\mathcal{T}$ by $\alpha_t(S)=e^{it}S$ (which is implemented by the unitary flow $Z$ defined by $Z_t\xi_k=e^{ikt}\xi_k$). Then $\alpha$ has a unique KMS state for all inverse temperature $\lambda\geq0$.

Let $P=1-SS^*$ and let $e_{k,\ell}=S^kP(S^*)^\ell$ for $k,\ell=0,1,2,\ldots$. Then $(e_{k,\ell})$ forms a family of matrix units and spans an ideal $\K$ of $\mathcal{T}$ as a closed linear subspace such that $\mathcal{T}/\K\cong C(\T)$ (and $e_{k,\ell}\xi_\ell=\xi_k$). And $\mathcal{T}$ is a closed linear span of $1,S^k,(S^*)^k,\ k\in\N$ and $\K$. Note that $\alpha_t(e_{k,\ell})=e^{i(k-\ell)}e_{k,\ell}$ and the induced flow on $C(\T)$ is given by translations. Then the KMS state $\omega_\lambda$ with $\lambda>0$ is given by
$$
\omega_\lambda(e_{k,k})=e^{-k\lambda}(1-e^{-\lambda})
$$
and $\omega_\lambda(e_{k,\ell})=0$ for $k\not=\ell$ and $\omega_0|\K=0$ (and $\omega_\lambda(S^k)=0=\omega_\lambda((S^*)^k)$ for $k=1,2,\ldots$).
\end{example}

The above $\alpha$ is not AI since it does not have a KMS state for $\lambda<0$, which follows from $1=\omega_\lambda(S^*S)=\omega_\lambda(S\alpha_{i\lambda}(S^*))=e^{\lambda}\omega_\lambda(SS^*)\leq e^{\lambda}$.  It also follows directly since the quotient flow on $C(\T)$ is not AI. But note that $\alpha|\K$ is AI.

One can derive $\omega_\lambda|\K$ as follows:
$$\omega_\lambda(e_{k,k})=\omega_\lambda(e_{k,\ell}e_{\ell,k})
=\omega_\lambda(e_{\ell,k}e^{-(k-\ell)\lambda}e_{k,\ell})=e^{-k\lambda+\ell\lambda}\omega_\lambda(e_{\ell,\ell}),
$$
 which implies that $e^{k\lambda}\omega_\lambda(e_{k,k})$ is independent of $k$.

The Toeplitz algebra $\mathcal{T}$ is not an MF C$^*$-algebra since it has a non-unitary isometry. The following example gives an example of a flow on an MF C$^*$-algebra.

\begin{example}
Let $S$ be as in the previous example and let $U=S^*\oplus S$ on $\Hil\oplus \Hil$. Let $A=C^*(U)$. Let $P_1=1-U^*U=P\oplus 0$ and $P_2=1-UU^*=0\oplus P$, where $P=1-SS^*$. Let $e^{(1)}_{k,\ell}=(U^*)^kP_1U^\ell=e_{k,\ell}\oplus 0$ and $e^{(2)}_{k,\ell}=U^kP_2(U^*)^\ell=0\oplus e_{k,\ell}$. The closed linear span of $e^{(i)}_{k,\ell}$ forms an ideal $\K_i$ for $i=1,2$ such that $\K_1\K_2=0$ and $A/(\K_1+\K_2)\cong C(\T)$. We define a flow on $A$ by $\alpha_t(U)=e^{it}U$ (which is implemented by the unitary flow $Z^*\oplus Z$). Then $\alpha$ has KMS states for all inverse temperature; if $\lambda>0$ the KMS state $\omega_\lambda$ satisfies $\|\omega_\lambda|K_2\|=1$ and if $\lambda<0$ then $\|\omega_\lambda|\K_1\|=1$ and $\omega_0|\K_1+\K_2=0$ for $\lambda=0$. Note that $\alpha$ is not approximately inner since the induced flow on the quotient $C(\T)$ is non-trivial.
\end{example}

To show that the above flow $\alpha$ is an MF flow on $A=C^*(U)$, let $B=\prod_n M_{n+1}$, $\beta_t=\prod_n \beta_{n,t}$, and $u=(u_n)\in B_\beta$ where $u_n=\sum_{k=1}^{n} e_{k+1,k}\in M_{n+1}$ and $\beta_{n,t}=\Ad(\sum_{k=1}^{n+1}e^{ikt}e_{k,k})$ on $M_{n+1}$. Then it follows that the C$^*$-subalgebra of $B_\beta/I$ generated by $u$ is isomorphic to $A$ with $\beta_t(u)=e^{it}u$. This follows because $p_1=1-u^*u=(e_{n+1,n+1})_n$ and $p_2=1-uu^*=(e_{1,1})_n$ are abelian projections in $B_\beta$ and $(u^*)^kp_1u^\ell=(e_{n+1-k,n+1-\ell})_n$ and $u^kp_2(u^*)^\ell=(e_{1+k,1+\ell})_n$ are families of matrix units mutually orthogonal in $B_\beta/I$.

%%%%%%%%%%%%%%%%%%%%%%%%%%%%%%%%%%%%%%%%%%%%%%%%%%%%%%%%%%%%%%%%%%%
\section{The case $F=M_k$}

We now turn to the examples described in the Introduction; translation-invariant flows on $A=\bigotimes_{m\in \Z} A(m)$ with $A(m)=F$ and their extensions to  $A\times_\gamma\Z$.

First we give a general result:

\begin{lem}
Let $A$ be a unital C$^*$-algebra. Let $\gamma$ be an automorphism of $A$ and $\alpha$ an AI flow on $A$ such that $\gamma\alpha_t=\alpha_t\gamma$. Suppose that there are a sequence $(h_n)$ in $A_{sa}$ and  central sequences $(u_n)$ and $(v_n)$ in $U(A)$ such that
\begin{enumerate}
\item  $\alpha_t(x)=\lim_n \Ad\,e^{ith_n}(x)$ uniformly in $t$ on a bounded  interval of $\R$ (or equivalently for each $t\in\R$);
\item  $\lim_n\|\Ad\,u_n\gamma(h_n)-h_n\|=0$;
\item  $\lim_n\|u_n-v_n\gamma(v_n)^*\|=0$.
\end{enumerate}
Then it follows that $\bar{\alpha}$ is AI on $A\times_\gamma\Z$. Here $\bar{\alpha}$ is defined by $\bar{\alpha}_t|A=\alpha_t$ and $\bar{\alpha}_t|C^*(\Z)=\id$, where $C^*(\Z)$ is the canonical C$^*$-subalgebra of $A\times_\gamma\Z$ coming from the action of $\Z$.
\end{lem}
\begin{pf}
Let $w_n=v_n^*u_n\gamma(v_n)$, which converges to $1$ as $n\ra\infty$. Then $\lim_n[v_n^*h_nv_n, w_nU]=0$ since $[\Ad\,v_n^*(h_n),w_nU]=\Ad\,v_n^*([h_n,\Ad\,v_n(w_nU)])$ and $\Ad\,v_n(w_nU)=v_nw_nUv_n^*=v_nw_n\gamma(v_n)^*U=u_nU$.

Denote by $\delta_\alpha$ the generator of $\alpha$. If $x\in D(\delta_\alpha)$ then there is a sequence $(x_n)$ in $A$ such that $\lim_nx_n=x$ and $\lim_n\ad(ih_n)(x_n)=\delta_\alpha(x)$ (coming from the fact that $\delta_\alpha$ is the graph limit of $\ad\,ih_n$, see, e.g., \cite{BR1}). Since $(v_n)$ is central we deduce that $\lim_n\Ad\,v_n^*(x_n)=x$ and $\lim_n i[\Ad\,v_n^*(h_n), \Ad\,v_n^*(x_n)]=\lim_n\Ad\,v_n^*(i[h_n, x_n])=\delta_\alpha(x)$, which implies that the graph limit of $\ad(iv_n^*h_nv_n)$ is $\delta_\alpha$ on $A$. Thus we conclude that the sequence $\ad(iv_n^*h_n v_n)$ converges to $\bar{\delta}_{\alpha}$ on the $^*$-algebra $D_0$ generated by $D(\delta_\alpha)$ and $U$ as a graph limit. Since $D_0$ is a core for $\bar{\delta}_\alpha$, this implies that the graph limit of $\ad(i v_n^*h_nv_n)$ is $\bar{\delta}_\alpha$; thus $\Ad\,e^{itv_n^*h_nv_n}(x)$ converges to $\bar{\alpha}_t(x)$ uniformly in $t$ on every bounded set of $\R$ for $x\in A\times_\gamma\Z$, i.e., $\bar{\alpha}$ is AI.
\end{pf}

We invoke the following two results on AI derivations, which may not be familiar with some readers now, to avoid giving the impression the above proof is a little sloppy. Here we call a derivation $\delta$ on a C$^*$-algebra $A$ an {\em AI derivation} if $\delta$ is the graph limit of $\ad\,ih_n$ for some sequence $(h_n)$ in $A_{sa}$. For a sequence in $A_{sa}$ to define an AI derivation we only need to check its possible 'domain' is dense; and for an AI derivation $\delta$ to generate a flow we only need to check whether the range of $\id\pm \delta$ is dense or not. These results are essentially found in Sakai's book \cite{Sak}.

\begin{prop}
Let $(h_n)$ be a sequence in $A_{sa}$ and define
$$
G=\{(x,y)\in A\times A\ |\ x_n\ra x,\ \ad ih_n(x_n)\ra y\}.
$$
Suppose that $D=\{x\ | \exists y\in A\ \ (x,y)\in G\}$ is dense in $A$. Then $G$ is the graph of a AI derivation $\delta$; $\D(\delta)=D$ and $G=\{(x,\delta(x))\ |\ x\in D\}$.
\end{prop}
\begin{pf}
Note that $G$ is a closed subspace of $A\times A$ such that $(x,y)\in G$ implies $(x^*,y^*)\in G$. If we have shown that $(0,a)\in G$ implies $a=0$ then $G$ is the graph of a closed linear map $\delta$ with $D(\delta)=D$. It is easy to show that $\delta$ is a derivation.

Note that $I=\{ y\in A\ |\ (0,y)\in G\}$ is a closed ideal of $A$. Because if $y\in I$ and $x\in D$ then it follows that $xy,yx\in I$. (If $a_n,x_n\in A$ satisfies that $a_n\ra0$, $\ad\,ih_n(a_n)\ra y$, $x_n\ra x$, and $\ad\,ih_n(x_n)\ra z$ for some $z\in A$, then $a_nx_n\ra 0$ and $\ad\,ih_n(a_nx_n)\ra yx$.) Since $D$ is dense and $I$ is a closed subspace, this implies that $I$ is an ideal.

If $I$ is non-zero there is a $y\in I_+$ such that $\|y\|=1$. Hence there is a sequence $(x_n)\in A$ such that $x_n\ra0$ and $\ad ih_n(x_n)\ra y$. We may suppose that $x_n=x_n^*$. Since $D$ is dense in $A$ there is, for any small positive $\epsilon<1/2$, $(a,b)\in G$ such that $a,b\in A_{sa}$ and $\|a-y\|<\epsilon$, which entails $\max\Spec(a)=\|a\|>1-\epsilon$.  Hence there is a sequence $(a_n)\in A_{sa}$ such that $a_n\ra a$ and $\ad ih_n(a_n)\ra b$. For any $\lambda\in\R$ it follows that $a_n+\lambda x_n\ra a$; thus there is a pure state $\phi_n$ such that $\phi_n(a_n+\lambda x_n)=\|a_n+\lambda x_n\|$ for sufficiently large $n$. Let $\phi_\lambda$ be an accumulation point of $(\phi_n)$. Then it follows that $\phi_\lambda(a)=\|a\|$ and $\phi_\lambda(b+\lambda y)=0$, which follows from $\phi_n\ad ih_n(a_n+\lambda x_n)=0$. Hence $0=\phi_\lambda(b+\lambda y)\leq \|b\|+\lambda \phi_\lambda(y)$ and $\phi_\lambda(y)=\phi_\lambda(a)+\phi_\lambda(y-a)\geq \|a\|-\|y-a\|>\|a\|-\epsilon>1-2\epsilon>0$, which is a contradiction for $\lambda<-\|b\|/(1-2\epsilon)$. See the proof of 3.2.9 of \cite{Sak}.
\end{pf}

\begin{prop}\label{well-behaved}
If $\delta$ is an AI derivation in $A$ then it is closed and satisfies
$$
\|x+\lambda\delta(x)\|\geq \|x\|,\ \ x\in\D(\delta),
$$
for all $\lambda\in\R$.
\end{prop}
\begin{pf}
A graph limit is closed if it is well-defined.
Suppose that $\delta$ is the graph limit of $\ad\,ih_n$ with $h_n\in A_{sa}$. Let $x\in \D(\delta)$ be a positive element. Then there is a sequence $(x_n)$  in $A_{sa}$ such that $x_n\ra x$ and $\ad\,ih_n(x_n)\ra \delta(x)$. Let $\phi_n$ be a state of $A$ such that $\phi_n(x_n)=\|x_n\|$ for large $n$ (which exists because $\|x_n\|=\max \Spec(x_n)$ for large $n$) and let $\phi$ be a weak$^*$ accumulation point, which automatically satisfies that $\phi(x)=\|x\|$. We may suppose that $\phi_n\ra \phi$; then $\phi(\delta(x))=\lim_n \phi_n(i[h_n,x_n])=0$ because $\phi_n$ is a character when restricted to $C^*(x_n)$, the C$^*$-algebra generated by $x_n$, i.e., for any positive $x\in\D(\delta)$ there is a state $\phi$ of $A$ such that $\phi(x)=\|x\|$ and $\phi(\delta(x))=0$. Let $x\in\D(\delta)$ and let $\phi$ be a state of $A$ such that $\phi(x^*x)=\|x\|^2$ and $\phi\delta(x^*x)=0$. Then if $\lambda\in\R$ then $\phi((x+\lambda\delta(x))^*(x+\lambda\delta(x)))=\|x\|^2 +\lambda\phi(\delta(x)^*x+x^*\delta(x))+\lambda^2\phi(\delta(x)^*\delta(x))\geq \|x\|^2$, i.e., $\|x+\lambda\delta(x)\|\geq \|x\|$. See 3.2.19 of \cite{Sak}.
\end{pf}

Now we consider the situation discussed at the end of Section 1 and we set $F=M_2$; $\gamma$ is the translation automorphism of $A=\bigotimes_{m\in\Z}A(m)$ with $A(m)=M_2$. We define $S\in M_2\otimes M_2$ by
$$
S=\sum_{i,j} e_{i,j}\otimes e_{j,i},
$$
where $e_{i,j},\ i,j=1,2$ are a family of matrix units of $M_2$. Note that $S$ is a self-adjoint unitary such that $\Ad\,S(x\otimes y)=y\otimes x$. Let $u_n$ denote $S$ in $A(-n)\otimes A(n+1)$.  Then $(u_n)$ is a central sequence of unitaries and $\Ad\,u_n$ switches $A(-n)$ and $A(n+1)$ and does nothing on the other factors. Thus $\Ad\,u_n\gamma$ cyclically permutes $A(-n),A(-n+1),\ldots,A(n)$ in particular. Let $\gamma^{(n)}$ denote the automorphism of the cyclic permutation of  factors of $A_n=\bigotimes_{m\in [-n,n]}A(n)$, i.e., $\gamma^{(n)}(A(m))=A(m+1)$ for $n\leq m<n$ and $\gamma^{(n)}(A(n))=A(-n)$. The system $(A_n,\gamma^{(n)})$ is based on $\Z_{2n+1}\equiv \Z/(2n+1)\Z$ in the same way $(A,\gamma)$ is based on $\Z$. We define a potential $\Phi_n$ over $\Z_{2n+1}$ from $\Phi$ over $\Z$ as follows: $\Phi_n(\Lambda')=(\gamma^{(n)})^m(\Lambda)$ if $\Lambda'=\Lambda+m\ {\rm modulo}\ 2n+1$ for some $\Lambda\subset [-n,,n]$ and $\Phi_n(\Lambda')=0$ otherwise. Then $\Phi_n$ is $\gamma^{(n)}$-invariant over $\Z_{2n+1}$. Let $h_n$ be the Hamiltonian over the periodic system $(A_n,\gamma^{(n)})$, i.e.,
$$
h_n=\sum_{\Lambda'}\Phi_n(\Lambda'),
$$
where the sum is taken over all $\Lambda'\subset \Z_{2n+1}$ and $A_n$ is identified with $A([-n,n])\subset A$. Then it follows that $\Ad\,e^{ith_n}$ converges to $\alpha^{\Phi}_t$ on $A$ and that $\gamma^{(n)}(h_n)=\Ad\,u_n\gamma(h_n)=h_n$. By the Rohlin property \cite{UHF-R} of $\gamma$ applied to $(u_n)$ we obtain a central sequence $(v_n)$ in $\U(A)$ such that $\|u_n-v_n\gamma(v_n)^*\|\ra0$. Hence the previous lemma implies that $\bar{\alpha}^{\Phi}$ is AI on $A\times_\gamma\Z$.

More generally one can prove the higher dimensional version. Let $d$ be a positive integer. Let $A(m)=M_k$ for each $m\in \Z^d$ and let $A=\bigotimes_{m\in \Z^d}A(m)$. We denote by $\gamma$ the action of $\Z^d$ by translations. If a $\gamma$-invariant potential $\Phi$ satisfies that $\sum e^{\lambda |\Lambda|}\|\Phi(\Lambda)\|<\infty$ for some $\lambda>0$ where the sum is taken over finite subsets $\Lambda\ni (0,\ldots,0)$, or $\Phi(\Lambda)$'s mutually commute, then one can define a flow $\alpha^\Phi$ on $A$ as before \cite{BR2}. One can also extend $\alpha^\Phi$ to a flow $\bar{\alpha}^\Phi$ on the crossed product $A\times_\gamma\Z^d$ imposing the condition $\bar{\alpha}_t^\Phi|C^*(\Z^d)=\id$.

\begin{theo}\label{UHF}
Let $A=\bigotimes_{m\in \Z^d} A(m)$ with $A(m)=M_k$ for some $k>1$ and $\gamma$ denote the action of $\Z^d$ by translations. Let $\alpha$ be the flow associated with a $\gamma$-invariant potential $\Phi$ with $\|\Phi\|_\lambda<\infty$ for some $\lambda>0$, or with $\Phi(\Lambda)$'s being commuting with each other. Then $\bar{\alpha}$ is AI on $A\times_\gamma\Z^d$.
\end{theo}
\begin{pf}
The case $d=1$ was treated in the above.

We concentrate on the case $d=2$; The case $d>2$ is similar. We also set $A(m)=M_2$ to spare the symbol $k$.

Let $n\in \N$. We identify $\Z_{2n+1}$ with $I_n=\{-n,-n+2,\ldots,n-1,n\}$ and let $\Phi_n$ be the potential on $\Z\times\Z_{2n+1}$ derived from $\Phi$ by imposing the {\em periodic boundary condition} on $I_n$. That is, for a finite subset $\Lambda$ of $\Z\times \Z_{2n+1}$ we define $\Phi_n(\Lambda)$ as $\gamma^{(n)}_{(0,k)}(\Phi(\Lambda'))$ if $\Lambda'\subset \Z\times I_n$ satisfies $\Lambda'+(0,k)=\Lambda$ (module $2n+1$ for the second coordinate) and $\Phi_n(\Lambda)=0$ if there is no such $\Lambda'$, where $\gamma^{(n)}$ is the action of $\Z\times\Z_{2n+1}$ on $A_n=\bigotimes_{m\in \Z\times \Z_{2n*1}}A(m)$. Let $\alpha^{(n)}$ denote the flow on $A_n$ associated with $\Phi_n$. By regarding $A_n=\bigotimes_{m\in \Z\times I_n}A(m)\subset A$ and $\alpha^{(n)}_t=\id$ on $\bigotimes_{m\not\in \Z\times I_n}A(m)\subset A$, we may suppose that $\alpha^{(n)}$ is a flow on $A$. Then it follows that $\alpha_t^{(n)}(x)$ converges to $\alpha_t^\Phi(x)=\alpha_t(x)$ uniformly in $t$ on every bounded set for all $x\in A$.

Regarding $A_n=\bigotimes_{k\in\Z}(\bigotimes_{\ell\in \Z_{2n+1}}A((k,\ell))$ we apply the argument for the case $d=1$. Let $\gamma_1^{(n)}$ (resp. $\gamma_2^{(n)}$) denote the restriction of $\gamma^{(n)}$ to $\Z\times \{0\}\cong\Z$ (resp. $\Z_{2n+1}$). Note that $\gamma_1^{(n)}$ has the Rohlin property even when restricted to $(A_n)^{\gamma_2}$, which is a unital simple AF C$^*$-algebra with unique tracial state. To prove the above theorem for $d=2$ we will need the following refinement:

\begin{lem}\label{periodic}
In the above situation $\bar{\alpha}^{(n)}$ is AI on $A_n\times_{\gamma^{(n)}_1}\Z$. Moreover there is a sequence $(h_m)$ in $(A_n^{\gamma^{(2)}})_{sa}$ and a sequence $(w_m)$ in $\U((A_n)^\gamma_n)$ such that $\bar{\alpha}^{(n)}_t(x)=\lim_m\Ad\,e^{ith_m}(x)$ uniformly in $t$ on every bounded set of $\R$ for $x\in A_n$, $w_m\ra 1$, and $[h_m,w_mU]=0$.
\end{lem}
\begin{pf}
Let $S=\sum_{i,j} e_{ij}\otimes e_{ji}$ in $M_2\otimes M_2$. Let $m$ be a large integer and let $S_k$ be the unitary in $A(-2m,k)\otimes A(2m+1,k)$ obtained as $S$ above. Let $u_m$ be the unitary $\bigotimes_{k\in Z_{2n+1}} S_k$ in $\bigotimes_{k\in Z_{2n+1}}(A(-2m,k)\otimes A(2m+1,k))$. Let $H_{m,n}$ denote the Hamiltonian on $\Z_{4m+1}\times \Z_{2n+1}$, i.e., $H_{m,n}=\sum_{\Lambda\subset \Z_{2m+1}\times \Z_{2n+1}}\Phi_{m,n}(\Lambda)$ where $\Phi_{m,n}$ is the potential obtained from $\Phi$ by imposing the periodic boundary condition on $I_m\times I_n$. Then $[H_m,u_mU_1]=0$, where $U_1$ is the canonical unitary in $A\times_{\gamma^{(n)}_1}\Z$. We note that $u_m$ and $H_{m,n}$ belong to $(A_n)^{\gamma_2}$.

Since the restriction of $\gamma^{(n)}_1$ to $(A_n)^\gamma$ has the Rohlin property (see \cite{AT} and \cite{LO05}) and $(u_m)$ is a central sequence in $A_n$, there is a sequence $(v_m)$ in $U((A_n)^\gamma)$ such that $\|u_m-v_m\gamma^{(n)}_1(v_m^*)\|\ra0$ and $(v_m)$ is central in $A_n$. (If we construct $v_m$ in a specific way as in \cite{AT} by using a path connecting $1$ to
$$
u_m\gamma(u_m)\gamma^2(u_m)\cdots \gamma^{\ell-1}(u_m)\in B\cap A(I_m\times \Z_{2n+1})'
$$
with $\ell=m,m-1$, then we would obtain $\|u_m-v_m\gamma_1^{(n)}(v_m^*)\|=O(1/m)$ and $v_m\in A(I_m\times \Z_{2n+1})'$.) Let $w_m=v_m^*u_m\gamma_{(1,0)}(v_m)\in B$, which implies that $w_mU_1=v_m^*u_mU_1v_m$. Then $\|w_m-1\|\ra0$ as $m\ra\infty$ and $v_m^*H_mv_m$ commutes with $w_mU_1$. Since $i[v_m^*H_mv_m,v_m^*xv_m]=i\Ad\,v_m^*[H_m,x]\ra \delta_\alpha(x)$ for $x\in A_{loc}=\bigcup_m A(I_m\times \Z_{2n+1})$, we may take $(v_m^*H_mv_m)$ for $(h_m)$ together with $(w_m)$.
\end{pf}

We identify $A_n$ as the C$^*$-subalgebra $A(\Z\times I_n)$ of $A=A(\Z^2)$. Let $S_i$ be the unitary in $A(i,-n)\otimes A(i,n+1)$ obtained as $S$ as in the proof of the above lemma. Fix $m\in\N$. We may assume that $h_m$ in Lemma \ref{periodic} is local, i.e., there is an $m'\in \N$ such that $h_m\in A(I_{m'}\otimes I_n)$. Let $\ell$ be an integer much greater than $m'$ and let $u_2$ denote the unitary $\bigotimes_{i=-\ell}^{i=\ell} S_i$. Since $h_m\in (A_n)^\gamma$ and $\Ad(u_2U_2)$ induces $\gamma^{(n)}_2$ on $A_n$ it follows that $[h_m, u_2U_2]=0$, where $U_2$ is the canonical unitary in $A\times_\gamma\Z^2$ corresponding to $\gamma_{(0,1)}$. Note that $u_2\gamma_{(1,0)}(u_2)^*=S_{-\ell} S_{\ell+1}^*\in A(\{-\ell,\ell+1\}\times \{-n,n+1\})$, which is central in $A$ as $\ell\ra\infty$. Since $\gamma_{(1,0)}$ has the Rohlin property, it follows that there is a unitary $v\in A$ such that $\|u_2\gamma_{(1,0)}(u_2)^*-v\gamma_{(1,0)}(v)^*\|$ is of the order $1/(\ell-m')$ and $v\in A(\{-m',-m'+1,\ldots,m'\}\times I_n)'$. Here we construct $v$ by using a path connecting $1$ to
$$
w\gamma_{(1,0)}(w)\gamma_{(2,0)}(w)\cdots \gamma_{(\ell-m'-1,0)}(w)=u_2\gamma_{(\ell-m',0)}(u_2)^*
$$
in $A([-\ell,-m']\times \{-n,n+1\})\cup [\ell+1,2\ell-m']\times\{-n,n+1\})$ etc. with $w=u_2\gamma_{(1,0)}(u_2)^*$. Then for $z=v^*u_2$ we obtain that $[h_m, zU_2]=[h_m,v^*]u_2 U_2+v^*[h_m,u_2U_2]=0$ and $\|z-\gamma_{(1,0)}(z)\|$ is of the order $1/(\ell-m')$.

For each $n\in\N$ we choose sequences $(h_m^{(n)})$ and $(w_m^{(n)})$ as in Lemma \ref{periodic}. Then we choose a subsequence $(m_n)$ such that $\ad\,ih_{m_n}^{(n)}$ converges to $\bar{\delta}_\alpha$ as a graph limit on $C^*(A,U_1)=A\times_{\gamma_1}\Z\subset A\times_\gamma\Z^2=C^*(A,U_1,U_2)$. In particular $w_{1,n}=w_{m_n}^{(n)}\in \U(A)$ satisfies  that $w_{1,n}\ra1$ and $\ad\,ih_{m_n}(w_{1,n} U_1)\ra0$. We will then choose $z_n$ in $U(A)$ such that $[h_{m_n}^{(n)},z_nU_2]=0$ and $\|z_n-\gamma_{(1,0)}(z_n)\|\ra0$ and $(z_n)$ is a central sequence in $A$. Now we denote $h_{m_n}^{(n)}$ by $h_n$.

Note that $\gamma$ has the Rohlin property (as an action of $\Z^2$; see \cite{Naka}). Hence there is a central sequence $(v_n)$ in $U(A)$ such that $\|z_n-v_n\gamma_{(0,1)}(v_n)^*\|\ra0$ and $\|v_n-\gamma_{(1,0)}(v_n)\|\ra0$. Let $w_{2,n}=v_n^*z_n\gamma_{(0,1)}(v_n)$, which converges to $1$ and satisfies that $[v_n^*h_nv_n, w_{2,n}U_2]=[v_n^*h_nv_n,v_n^*z_nU_2v_n]=\Ad(v_n^*)[h_n,z_nU_2]=0$. Note also that $[v_n^*h_nv_n,v_n^*w_{n,1}U_1v_n]\ra0$ and $v_n^*w_{n,1}U_1v_n= v_n^*w_{n,1}\gamma_{(1,0)}(v_n)U_1 \ra U_1$ (i.e., the choice of $v_n$ for $z_n$ is done to preserve the property of near $\gamma_{(1,0)}$-invariance for $z_m$). Thus one can conclude that $\ad\,iv_n^*h_nv_n$ converges to $\bar{\delta}_\alpha$ as a graph limit on $A\times_\gamma\Z^2$.
\end{pf}

We will give another proof to the above theorem which seems to work only for the case $d=1$, but which can be extended to $F$ with non-trivial center. What we did to show that $\bar{\alpha}$ is AI on $A\times_\gamma\Z$ was to find a sequence $(h_n)$ in $A$ and a sequence $(w_n)$ in $\U(A)$ such that $\Ad\,e^{ith_n}\ra \alpha_t$ on $A$, $\|w_n-1\|\ra0$,  and $\|\Ad\,w_n\gamma(h_n)-h_n\|\ra0$. (Whether this is necessary or not is unknown except for the case $F$ is a full matrix algebra \cite{AI-MF}.) What we will do below is to find, for any $\epsilon>0$, some $N\in\N$ and sequences $h_n=(h_{n,0},h_{n,1},\ldots,h_{n,N-1})$ with $h_{n,i}\in A_{sa}$ and $w_n=(w_{n,0},w_{n,1},\ldots,w_{n,N-1})$ with $w_{n,i}\in \U(A)$ such that $\Ad\,e^{ith_{n,i}}\ra\alpha_t$ on $A$, $\|w_{n,i}-1\|<\epsilon$, and $\|\Ad\,w_{n,i}\gamma(h_{n,i})-h_{n,i+1}\|<\epsilon$ with $h_{m,N}=h_{m,0}$. Here $N$ should be such that there is a sequence of projections $e_{n,0},e_{n,1},\ldots e_{n,N-1}$ in $A$  such that $\sum_i e_{n,i}=1$, $\|\gamma(e_{n,i})- e_{n,i+1}\|\ra0$ with $e_N=e_0$, and $(e_{n,i})$'s are central sequences. Then by setting $h_n=\sum_i e_{m,i}h_{n,i}e_{m,i}$ and $w_n$ the unitary part of $\sum_i w_{n,i}\gamma(e_{n,i})\approx 1$ for sufficiently large $n$ and then large $m$ we would obtain that $\Ad\,e^{ith_n}\approx\alpha_t$, $\|w_n-1\|<\epsilon$, and    $\|\Ad\,w_n\gamma(h_n)-h_n\|<\epsilon$; thus we will be reduced to the first case.

We first assume that $\Phi$ is a finite-range potential, i.e., there is a constant $r$ such that if $\Lambda\subset \Z$ has two points $k,\ell$ such that $|k-\ell|>r$ then $\Phi(\Lambda)=0$; the minimum of such $r$ will be called the {\em range} of $\Phi$.

Let $k_0\in \N$ such that $N=2^{k_0}\gg 2r$. For $k\not=\ell$ let $S(k,\ell)$ denote the self-adjoint unitary $S$ in $A(k)\otimes A(\ell)$; $\Ad\,S(k,\ell)$ switches $A(k)$ and $A(\ell)$. Let
$$
V=\prod_{i=1}^{2N} S(-7N+i,N+i),
$$
which is a self-adjoint element whose adjoint action switches $A[-7N+1,-5N]$ and $A[N+1,3N]$. Let $E$ be the projection with $V=e^{i\pi E}$ and set $u=e^{i\pi E/2N}\in A[-7N+1,3N]$, which satisfies $\|u-1\|\leq \pi/2N$.

For $k\in\Z$ let
$$
W(k)=\sum_\Lambda \Phi(\Lambda),
$$
where the summation is taken over all $\Lambda$ with $\Lambda\cap (-\infty,k]\not=\emptyset$ and $\Lambda\cap [k+1,\infty)\not=\emptyset$, i.e., $W(k)$ denotes the {\em interaction} between the two regions $(-\infty,k]$ and $[k+1,\infty)$. We set $w(k)=W(k)/2N$. Note that $\gamma(w(k))=w(k+1)$.

We set
$$
h_0=H[-5N+1,-3N]+H[-3N+1,3N]=H[-5N+1,3N]-W(-3N).
$$
Then we define
\begin{align*}
h_1&=\gamma(u)\gamma\big(h_0+w(-3N)-w(N)\big)\gamma(u)^*\\
   &=\gamma(u)\big(\gamma(h_0)+w(-3N+1)-w(N+1)\big)\gamma(u)^*
\end{align*}
which belongs to $A[-7N+2,3N+1]$ and satisfies that $\|h_1-\gamma(u)\gamma(h_0)\gamma(u)^*\|\leq 2\|W(0)\|/2N$. We further define
\begin{align*}
h_2 &=\gamma^2(u)\gamma\big(h_1+\gamma(u)(w(-3N+1)-w(N+1))\gamma(u)^*\big)\gamma^2(u)^*\\
    &=\gamma^2(u)^2\big(\gamma^2(h_0)+2w(-3N+2)-2w(N+2)\big)\gamma^2(u)^{-2}.
\end{align*}
Then it follows that $h_2\in A[-7N+3, 3N+2]$ and $\|h_2-\gamma^2(u)\gamma(h_1)\gamma^2(u)^*\|\leq 2\|W(0)\|/2N$. Similarly we define $h_k,\ k=3,4,\ldots,2N$ by the relation
$$
h_k-\Ad\gamma^k(u)(\gamma(h_{k-1}))=\Ad\gamma^k(u)^k(w(-3N+k)-w(N+k))
$$
which amounts to
$$
h_k=\gamma^k(u)^k\big( \gamma^k(h_0)+kw(-3N+k)-kw(N+k)\big) \gamma^k(u)^{-k}.
$$
Then $h_{2N}$ equals to $h_0$ because
\begin{equation*}
\begin{array}{ll}
 h_{2N}&= \gamma^{2N}(u)^{2N}\big( \gamma^{2N}(h_0)+W(-N)-W(3N) \big)\gamma^{2N}(u)^{-2N}\\
 &=\gamma^{2N}(u)^{2N}\big(H[ -3N+1,5N]-W(3N)\big)\gamma^{2N}(u)^{-2N}\\
 &=\gamma^{2N}(V)\big(H[ -3N+1,3N]+H[3N+1,5N]\big)\gamma^{2N}(V)^*\\
 &=H[-5N+1,-3N]+H[-3N+1,3N].
 \end{array}
\end{equation*}
Let $\epsilon>0$ and note that $N=2^{k_0}$. By using the Rohlin property of $\gamma$, let $e_0,e_1,\ldots,e_{2N-1}$ be a family of projections in $A$ such that $\sum_i e_i=1$, $e_i\in A[-7N,5N]'$, and $\|\gamma(e_i)-e_{i+1}\|<\epsilon$ with $e_{2N}=e_0$. We set
$$
h=\sum_{i=0}^{2N-1} h_ie_i
$$
and
$$
v=\sum_{i=0}^{2N-1} \gamma^{i+1}(u)\gamma(e_i).
$$
Note that $\|v-1\|\leq \pi/2N$. Since $h_i,\gamma(h_i),\gamma^i(u)\in A[-7N,5N]$ it follows that
\begin{equation*}
\begin{split}
\|h-v\gamma(h)v^*\|&\leq \sum_{i=0}^{2N-1}\|h_{i+1}e_{i+1}-h_{i+1}\gamma(e_i)\|+
  \|\sum_{i=0}^{2N-1} (h_{i+1}-\gamma^{i+1}(u)\gamma(h_i)\gamma^{i+1}(u)^*)\gamma(e_i)\|\\
  &\leq (\sum_i \|h_{i+1}\|)\epsilon+\|W(0)\|/N.
\end{split}
\end{equation*}
By choosing $N$ sufficiently large and then choosing $\epsilon$ sufficiently small one would obtain that $\|[h,vU]\|=\|h-v\gamma(h)v^*\|\approx0$.
Since $\ad(ih)$ equals $\delta_{\alpha^\Phi}$ on $A[-N+1+r,N-r]$, we can obtain the conclusion of Theorem \ref{UHF} for a finite-range $\Phi$ in the case $d=1$.

In general we approximate $\Phi$ by finite-range $\Phi_n$ defined as
$\Phi_n(\Lambda)=\Phi(\Lambda)$ if ${\rm dim}(\Lambda)\leq n$ and $=0$ otherwise, where ${\rm dim}(\Lambda)=\max \{|k-\ell|\ ; \ k,\ell\in \Lambda\}$. Then we reach the conclusion since $\bar{\alpha}^{\Phi_n}_t\ra \bar{\alpha}_t^\Phi$

%%%%%%%%%%%%%%%%%%%%%%%%%%%%%%%%%%%%%%%%%%%%%%%%%%%%%%%%%%%%%%%%%%%%%%%%%%%%%%%%%%%%%%%%%%

\section{The case $F=M_k\oplus M_\ell$}

As a typical non-factor case we will consider the case $F=M_k\oplus M_\ell$; the arguments below will be similar if $F$ contains more than two factors. Since the values of $k$ and $\ell$ do not affect the arguments below up to the last point we will assume $F=M_2\oplus M_3$ to free the symbols $k,\ell$ for other uses.

The center of $A=\bigotimes_{n\in \Z}A(n)$ with $A(n)=F$ is  non-trivial, which we naturally identify with the continuous functions on $X=\prod_{n\in \Z}\{2,3\}$. We denote by $\gamma$ the translation automorphism of $A$.

Suppose that we are given a $\gamma$-invariant finite-range potential $\Phi$ for $A$; we denote by $\alpha^\Phi$ the flow on $A$ induced from $\Phi$ and by $\delta^\Phi$ the generator of $\alpha^\Phi$. Let $r$ be the range of $\Phi$. Then for it follows that for $x\in A[-N,N]$
$$
\delta^\Phi(x)=i\sum_{\Lambda\subset [-N-r,N+r]} [\Phi(\Lambda),x]
$$
and $\delta^\Phi$ is the closure of the restriction of $\delta^\Phi$ to $A_{loc}=\bigcup_N A[-N,N]$.

Let $N\in\N$ be  such that $N>r$. We will identify $X$ with $\prod_{i\in\Z} Y_i$ with $Y_i=\{2,3\}^{2N}=\prod_{-N< j\leq N}\{2,3\}$ by sending $x\in X$ to $y=(y_i)$ with $y_i=(x_{2Ni+j})_{-N< j\leq N}$. Let $M\in\N$ be much larger than $2^{2N}$ which is the cardinality of $Y_i$. (If $F$ contains three factors then $M$ should be much larger than $3^{2N}$.)

Fix $x=(x_k)_{-M\leq k\leq M}\in \prod_{-M\leq k\leq M} Y_k$. For each $y\in Y\equiv\{2,3\}^{2N}$ let
$$
I_y=\{k\ |\ -M\leq k\leq M,\ x_k=y\}.
$$
When $I_y$ is non-empty we will construct a directed graph whose vertices are elements of $I_y$ and whose directed edges or arrows are specified only by sources and targets. Let $I_0=I_y$ and $J_0=I_y$ and we draw an arrow from an element of $I_0$ to another of $J_0$ in the following way. Take a pair of $k\in I_0$ and $\ell\in J_0$ such that
$$
4\leq k-\ell\leq M.
$$
If there are more than one such pairs we choose one with $k$ maximal and then $\ell$ maximal and draw an arrow from $k$ to $\ell$. If none, we stop here, i.e., our directed graph has only vertices $I_y$ and no arrows. We then apply the same procedure to $I_1=I_y\setminus \{k\}$ and $J_1=I_y\setminus\{\ell\}$, i.e., if there are more than one pairs $k\in I_1$ and $\ell\in J_1$ such that $4\leq k-\ell\leq M$ then we choose one with maximal $k$ and then with maximal $\ell$ and draw an arrow from $k$ to $\ell$. If none we stop here. If we have drawn an arrow we apply the same procedure to $I_2=I_1\setminus \{k\}$ and $J_2=J_1\setminus\{\ell\}$ until there are no such pairs $k\in I_i$ and $\ell\in J_i$ for some $i$; if it stops after defining $I_i$ and $J_i$ we have drawn $i$ arrows among the points in $I_y$. Each point in $I_y$ has at most one arrow leaving it and at most one arrow targeting it.

We claim that $B= \{k\in I_y\ |\ k>0,\ k\ {\rm\ is\  not\ a\ source}\}\subset \{1,2,\ldots,M\}$ has at most four points. Suppose that $B$ has more than four points and let $k_{i},\ i=1,2,\ldots,m$ be the enumeration of $B$ in descending order with $m>4$. Since $M\geq k_{1}-k_{5}\geq 4$, the reason why we do not have an arrow from $k_{1}$ to $k_{5}$ must be that $k_5$ is already the target of an arrow starting from some $k_0>k_{1}$. Then since $k_0>k_{1}>k_{2}>k_{3}>k_{4}$, the inequality $k_0-k_{4}\geq 4$ shows that $k_{4}$ is the target of an arrow leaving $k_{-1}>k_0$. We continue this way ad infinitum; if $k_{n-5}$ is introduced as the source of the arrow to $k_n$, the reason why no arrow exists from $k_{n-5}$ to $k_{n-1}$ must be that $k_{n-1}$ is already the target from an arrow leaving from $k_{n-6}>k_{n-5}$. This is absurd since we are dealing with a subset of integers between $1$ and $M$.

We also claim that $B'=\{k\in I_y\ |\ k<0, \ k {\rm\ is\ not\ a\ target}\}$ has at most four points. Suppose that $B'$ has more than four points and let $k_i,\ i=1,2,\ldots,m$ be the enumeration of $B'$ in descending order with $m>4$. Since $k_1-k_5\geq 4$ the reason $k_5$ is not the target of an arrow from $k_1$ is that the arrow from $k_1$ has landed on $\ell$ with $k_5<\ell<k_1$ and $\ell\not=k_{2},k_3,k_4$. Let $S_1=\{\ell,k_2,k_3,k_4\}$. Then the maximum $\ell_1$ of $S_1$ is greater than $k_5$ at least by 4. The reason why $k_5$ is not the target of an arrow from $\ell_1$ is that the arrow from $\ell_1$ must have landed on a value greater than $k_5$, which is none of $S_1\setminus \{\ell_1\}$ (since $\ell$ is already a target and $k_2,k_3,k_4$ are not targets). Let $S_2$ denote the union of $S_1\setminus \{\ell_1\}$ and the singleton consisting of the target of $\ell_1$. Then the maximum $\ell_2$ of $S_2$ is greater than $k_5$ by at least 4. Since no elements of $S_2$ can be a target of $\ell_2$ we again obtain the target of an arrow from $\ell_2$, which is greater than $k_5$. By adding this element to $S_2\setminus \{\ell_2\}$ we define $S_3$ and continue this argument ad infinitum to reach a contradiction.

Note that $[-M,M]$ is the union $I_y$ over $y\in \prod_{N<k\leq N}\{2,3\}$ and we take the union of the directed graphs based on $I_y$ over all $y\in Y$.

We denote by $e(x)$ the central projection of $A$ corresponding to
$$
x=(x_k)_{-M\leq k\leq M}\in \prod_{-M\leq k\leq M}Y_k=\prod_{-(2M+1)N<i\leq (2M+1)N}\{2,3\}
$$
being regarded as a cylinder set of $X$. We specify
$$
h_0(x)\in \big(\bigotimes_{-(2M+1)N<n\leq (2M+1)N}A(n)\big) e(x)\subset Ae(x)
$$
by first giving $h(k)\in \bigotimes_{-N<j\leq N}A(2Nk+j),\ k=-M,-M+1,\ldots, 0,\ldots M$ as follows:

If $k<-2$ is a target, then set $h(k)=0$; if $k<-2$ is not a target then
$$h(k)=\frac{M+k}{M-2}H[(2k-1)N+1,(2k+1)N];
$$
if $k=-2,-1,0,1$ set $h(k)=H[(2k-1)N+1,(2k+1)N]$; if $k\geq 2$ is a source and if its target is greater than or equal to $-2$, set $h(k)=0$; if $k\geq 2$ is a source but if its target is smaller than  $-2$, set $h(k)=H[(2k-1)N+1,(2k+1)N]$; if $k\geq 2$ is not a source then
$$h(k)=\frac{M-k}{M-2}H[(2k-1)N+1,(2k+1)N].
$$
With those $h(k)$ we define
$$
h_0(x)=\big(\sum_{k=-M}^{M}h(k)+W(-N)+W(N)\big)e(x),
$$
which has the {\em right} Hamiltonian in $[-3N+1,3N]$. (The idea is the following: If $k<-2$ is a target  the Hamiltonian, $h(k)$, on $[(2k-1)N+1,(2k+1)N]$ will be inherited from the source; otherwise we have to build it up a bit by a bit as $k$ increases. This inheritance occurs when $k=-2$. For $k=-2,-1,0,1$ we assign  the {\em right} Hamiltonian to $h(k)$; If $k\geq 2$ is not a source we have to diminish $h(k)$ a bit by a bit as $k$ increases. Otherwise $h(k)$ will be transferred to the target when it is $-2$.)

Note that the number of $\{k\, |2\leq k\leq M, \ h(k)\not=0\}$ (resp. $\{k\,|-M\leq k<-2, \ h(k)\not=0\}$) is at most $4\cdot 2^{2N}$ and $h_0(x)$ depends only on $x_k,\ -(M-1)\leq k\leq M-1$ as $h(-M)=0=h(M)$ and $h(k),\ -M<k<-2$ depends on $x_\ell,\ -M<k+4\leq \ell\leq k+M<M$ and $h(k),\ 2\leq k<M$ depends on $-M<k-M\leq \ell\leq k-4<M$.

We denote by $E_-(x)$ the sum of $H[(2k-1)N+1,(2k+1)N]$ over the $k<-2$ which are not targets. We also denote by $E_+(x)$ the sum of $H[(2k-1)N+1,(2k+1)N]$ over the $k>2$ which are not sources. Note that $\|E_\pm (x)\|\leq 4\cdot 2^{2N}\|H[-N+1, N]\|$. Letting $C=\sum_{X\ni 0}\|\Phi(X)\|\leq \|\Phi\|_\lambda<\infty$ one has the estimate $\|H[-N+1,N]\|\leq 2NC$. Thus it follows that
$$
\|E_{\pm}(x)\|\leq 8\cdot 2^{2N}NC.
$$

%As far as  $y\in X$ satisfies $e(y)\gamma^{2N}(e(x))\not=0$, $h_0(y)$ is independent of the choice of $y$. Thus $h_0(\gamma^{2N}(x))=\sum_y h_0(y)\gamma^{2N}(e(x))$ is also well-defined, where $\gamma^{2N}(x)$ is the $2N$-shift of $x$ to the right. Hence it follows that $\sum_xh_0(x)=\sum_y h_0(\gamma^{2N}(y))$.

If there is an arrow from $k_0\in [1,M-3]$ to $-3$ then let
$$V=\prod_{i=1}^{2N}S(-7N+i,(2k-1)N+i),
$$
which is a self-adjoint unitary whose adjoint action switches $A(-7N+i)$ and $A((2k_0-1)N+i)$ for $i=1,2,\ldots,2N$ simultaneously. Let $E$ be the projection with $V=e^{i\pi E}$ and set $u=e^{i\pi E/2N}\in A[-7N+1, (2k+1)N]$, which satisfies that $\|u-1\|\leq \pi/2N$. If there is no arrow to $-3$ then we set $u=1$. We denote by $\gamma$ the right shift homeomorphism of $X$ and by $$\gamma(x)\in \prod_{-(2M+1)N+1<k\leq (2M+1)N+1}\{2,3\}
$$
the right shift of $x\in \prod_{-(2M+1)N<k\leq (2M+1)N}\{2,3\}$. We define
$$h_1(\gamma(x))\in \big( \bigotimes_{-(2M+1)N+1<n\leq (2M+1)N+1}A(n)\big)\gamma(e(x))
$$
as follows:
\begin{align*}
& h_1(\gamma(x))- \Ad\gamma(u)\gamma(h_0(x)) \\
&=\Ad\gamma(u)\big( \frac{1}{2N(M-2)}\gamma(E_-(x)-E_+(x))  +\frac{1}{2N}(W(-3N+1)-W(N+1))\big)\gamma(e(x)).
\end{align*}
Then one computes:
$$
\|h_1(\gamma(x))-\Ad\gamma(u)\gamma(h_0(x))\| \leq \frac{8\cdot 2^{2N} C}{M-2} +\frac{\|W(0)\|}{N}.
$$
Thus for a suitable choice of $N$ and $M$ the right hand side can be made arbitrarily small. In general we define $$h_k(\gamma^k(x))\in \big( \bigotimes_{-(2M+1)N+k<n\leq (2M+1)N+k} A(n) \big)\gamma^k(e(x))
$$
for $k\in [1,2N]$ such that
\begin{align*}
& h_k(\gamma^k(x))-\Ad\gamma^k(u)\gamma(h_{k-1}(\gamma^{k-1}(x))) \\
&=
\Ad\gamma^k(u)^k \big( \frac{1}{2N(M-2)} \gamma^k(E_-(x)-E_+(x)) +\frac{1}{2N}(W(-3N+k)-W(N+k))\big) \gamma^k(e(x)).
\end{align*}
Then $h_k(\gamma^k(x))-\Ad\gamma^k(u)^k \gamma^k(h_0(x))$ equals
$$
\Ad\gamma^k(u)^k\big( \frac{k}{2N(M-2)}\gamma^k(E_-(x)-E_+(x)) +\frac{k}{2N}(W(-3N+k)-W(N+k))\big) \gamma^k(e(x)).
$$

Eventually $h_{2N}(\gamma^{2N}(x))$ equals
$$
\Ad\gamma^{2N}(u)^{2N}\gamma^{2N}\big(h_0(x) +\frac{1}{M-2}(E_-(x)-E_+(x))+W(-3N)-W(N)\big)\gamma^{2N}(e(x)),
$$
which is $\sum h_0(y)\gamma^{2N}(e(x))$ with the sum taken over $y\in \prod _{-(2M+1)N<i\leq (2M+1)N}\{2,3\}$ such that $y_i=x_{i-2N}, \ -(2M+1)N+2N<i \leq (2M+1)N$ or $e(y)\gamma^{2N}(e(x))\not=0$.

The reason is as follows: If $-3$ is not a target then $u=1$ and $h(-3)=\frac{M-3}{M-2}H[-7N+1,-5N]$. Hence we deduce that
\begin{align*}
& \gamma^{2N}\big(\sum_{k=-M}^{M} h(k)+W(-N)+W(N)+\frac{1}{M-2}(E_-(x)-E_+(x))+W(-3N)-W(N)\big)\\
&= \sum_{k=-M+1}^M h'(k)+W(-N)+W(N),
\end{align*}
where $h'(k)$ belongs to $A[(2k-1)N+1,(2k+1)N]$ and is the $h(k)$ for $\gamma^{2N}(x)$ in place of $x$. Note that $h'(k)$ appears as the $2N$-shift of a modification of $h(k-1)$, which is accommodated by $(E_-(x)-E_+(x))/(M-2)$. If $-3$ is a target, then there is an arrow from $k_0 (\leq M-3)$ to $-3$, and $\gamma^{2N}(u)^{2N}=\gamma^{2N}(V)$ whose adjoint action switches $A[-5N+1,-3N]$ and $A[(2k_0+1)N+1,(2k_0+3)N]$. Since $h(k_0)=H[(2k_0-1)N+1,(2k_0+1)N]$ and $h(-3)=0$ we deduce that $h'(-2)=H[-3N+1,-N]$, which is inherited from $\gamma^{2N}(h(k_0))$. The other terms are dealt with in the same way as above.

Since $h_{2N}(\gamma^{2N}(x))=\sum h_0(y)\gamma^{2N}(e(x))$ over all  $y\in \prod_{-(2M+1)N<i\leq (2M+1)N}\{2,3\}$ it follows that
$$
\sum_x h_{2N}(\gamma^{2N}(x))=\sum_x h_0(x).
$$
We shall formulate the above arguments as follows.

\begin{prop}\label{finite-d}
Let $F$ be a finite-dimensional C$^*$-algebra and let $A=\bigotimes_{n\in\Z}A(n)$ where $A(n)=F$. Let $\gamma$ denote the right-translation automorphism of $A$ and let $\Phi$ be a $\gamma$-invariant finite-range potential for $A$. Then for any $\epsilon>0$ and $K\in\N$ there is an $N_0\in\N$ satisfying: For any $N\geq N_0$ there are $h_k\in (A_{loc})_{sa}$ and $u_k\in \U(A_{loc})$ for $k=0,1,2,\ldots,N-1$ such that $\|h_{k+1}-\Ad\,u_k\gamma(h_k)\|<\epsilon$ with $h_N=h_0$, $\|u_k-1\|<\epsilon$, and $\Ad(ih_k)=\delta^\Phi$ on $A([-K,K])$ for all $k=0,1,\ldots,N-1$.
\end{prop}

In the above argument with $F=M_2\oplus M_3$ we set $h_i=\sum h_i(x)$ where the summation is over all $x\in \prod_{-(2M+1)N<k\leq (2M+1)N}\{2,3\}$ and $u_i=\gamma^{i+1}(u)$ for $i=0,1,2,\ldots,2N-1$. Then if $N$ and $M$ are sufficiently large, the families $h_i,u_i,\ i=0,1,2,\ldots,2N-1$ satisfy the condition $\|h_{i+1}-\Ad\,u_i\gamma(h_i)\|<\epsilon$ for $i=0,1,2,\ldots,2N-1$ and the condition $\Ad(ih_i)=\delta^\Phi$ on $A[-N+r+1,N-r])$ where $r$ is the range of $\Phi$. This argument can be generalized to cover an arbitrary $F$.

In our previous arguments there are $2N$ of $h_i$'s, which comes from choosing $[-N+1,N]$ as a basic interval, which has $2N$ points. We could choose as the basic interval an arbitrary interval of $\Z$ of any length if it contains $[-N_0+1,N_0]$ for a sufficiently large $N_0$.

\begin{theo}
Suppose that $F\cong M_k\otimes F_1$ with $k>1$ and $F_1$ finite-dimensional in the situation of Proposition \ref{finite-d}. Let $\Phi$ be a $\gamma$-invariant potential such that $\|\Phi\|_\lambda<\infty$ or $\Phi(\Lambda)$'s commute with each other. Then $\bar{\alpha}^\Phi$ is AI on $A\times_\gamma\Z$.
\end{theo}

In this case $\gamma$ has the {\em cyclic} Rohlin property, i.e., for any $n\in \N$ and $\epsilon>0$ there is a sequence $e_0,e_1, \ldots,e_{k^n-1}$ of projections in $A$ such that $\sum_ie_i=1$ and $\|\gamma(e_i)-e_{i+1}\|<\epsilon$ for $i=0,1,\ldots,k^n-1$ with $e_{k^n}=e_0$. We may suppose that $e_i$'s are local, i.e., $e_i\in A[-K,K]$ for some $K$. By replacing $e_i$ by $\gamma^L(e_i)$ with $L$ large we can suppose that $e_i$'s commute with any finite number of local elements.

By Proposition \ref{finite-d} for any $\epsilon>0$ and $K\in\N$ one can choose $h_i,u_i,\ i=0,1,2,\ldots,k^n-1$ as specified there. We choose projections $e_i,\ i=0,1,2,\ldots,k^n-1$ as in the above paragraph such that $e_i,\gamma(e_i)$ commute with $h_j,u_j$ for all $i,j$. Then we set
$$
h=\sum_i h_i e_i,\ \ \ u=\sum_i u_i\gamma(e_i).
$$
Then $h$ is self-adjoint, $\|u-1\|$ is small, $\ad(ih)=\delta_\alpha$ on $A[-K,K]$, and $\|h-\Ad\,u\gamma(h)\|\leq \sum_i\|h_{i+1}\|\|e_{i+1}-\gamma(e_i)\|+\max_i\|h_{i+1}-\Ad\,u_i\gamma(h_i)\|$ which can be assumed arbitrarily small.

\begin{cor}
Let $F$ be a simple AF C$^*$-algebra such that $[1]\in K_0(A)$ is divided by an integer $k>1$ and let $A=\bigotimes_{m\in\Z}A(m)$ with $A(m)=F$. Let $\Phi$ be a $\gamma$-invariant potential in $A$ such that $\|\Phi\|_\lambda<\infty$ for some $\lambda>0$ or $\Phi(\Lambda)$'s commute with each other, which generates a flow $\alpha^\Phi$ on $A$. Then $\bar{\alpha}^\Phi$ is AI on $A\times_\gamma\Z$.
\end{cor}

%%%%%%%%%%%%%%%%%%%%%%%%%%%%%%%%%%%%%%%%%%%%%%%%%%%%%%%%%%%%%%%%%%%%
\section{The case of no interactions}

The arguments in the previous section are just good enough to show that the flow $\bar{\alpha}$ on $A\times_\gamma\Z$ induced by $\alpha$ on $A$ is AI in the case $A=\bigotimes_\Z F$ with $F=M_k$ or $M_k\otimes F_1$ for $k\geq 2$. The flows on $A$ defined through potentials are actually continuously AI and the question remains of whether $\bar{\alpha}$ is also continuously AI. We shall show in this section the flows defined through potentials of no interaction (between different lattice points) induce continuously AI flows on $A\times_\gamma\Z$; the proof is a bit more elaborate but based on the similar idea as before.

We say the potential $\Phi$ on $A=\bigotimes_\Z F$ has no interactions if $\Phi(X)=0$ for $X$ containing more than one elements. In this case the flow $\alpha$ is defined as $\alpha_t=\bigotimes_\Z \Ad\,e^{ith}$ for some $h=h^*\in F$. To make things look a bit more non-trivial we will later assume that $F$ is a UHF algebra and $A(m)=F$ for all $m\in \Z$ and set $A=\bigotimes_{m\in \Z}A(m)$, which is again a UHF algebra. We denote by $\gamma$ the action of $\Z$ on $A$ by translations as before and we let $\beta$ be an AI flow on $F$ and define a flow $\alpha$ on $A$ by $\bigotimes_{\Z}\beta$. Our aim is to show that if $\beta$ is continuously AI then so is the extension $\bar{\alpha}$ of $\alpha$ to a flow on $A\times_\gamma\Z$.

First of all we prove the following:

\begin{prop}\label{no-interaction}
Let $A=\bigotimes_{m\in\Z}A(m)$ with $A(m)=M_k$ for some $k\geq 2$ and let $\gamma$ be the action of $\Z$ on $A$ by translations. Let $h\in M_k$ such that $h=h^*$ and $h_m\in A(m)$ be a copy of $h$ and define a flow $\alpha$ on $A$ by $\alpha_t=\lim_{n\ra\infty} \Ad \exp it\{\sum_{m=-n}^n h_m$\}. Then the extension $\bar{\alpha}$ to a flow on the crossed product $A\times_\gamma\Z$ is continuously AI, where $\bar{\alpha}_t|A=\alpha_t$ and $\bar{\alpha}_t(U)=U$ with $U$ the $\gamma$-implementing unitary of $A\times_\gamma\Z$.
\end{prop}

Define a unitary $S\in M_k\otimes M_k$ by $S=\sum_{i,j} e_{ij}\otimes e_{ji}$ where $e_{ij}, \, i,j=1,2,\ldots,k$ are a family of matrix units. Since $S$ is a self-adjoint unitary the determinant of $S$ is $1$ or $-1$. (Actually $\det(S)=(-1)^{k(k-1)/2}$.) We denote by $S(m,n)$ the $S$ sitting in $A(m)\otimes A(n)$.

We will choose an increasing sequence $(N_i)$ of positive integers. For each $i$ we define
\begin{multline*}
u_i=  S(N_{i+1}+1,N_{i+1}) S(N_{i+1},N_{i+1}-1)\cdots S(N_i+2,N_i+1)S(N_i+1,-N_i)\cdot \\
     S(-N_i,-N_i-1)      \cdots S(-N_{i+1}+1,-N_{i+1})
\end{multline*}		
in $A[-N_{i+1},-N_i]\otimes A[N_i+1, N_{i+1}+1]\subset A$. Since $u_i$ is a product of an odd number of $S$, if $k$ is odd then $\det(u_i)=\det(S)$ in this matrix subalgebra. (If $k$ is even then $\det(u_i)=1$. Note that $u_i$ is regarded as a matrix of degree $k^\ell$ for $\ell=2(N_{i+1}-N_i+1)$. )
Then $\Ad\,u_i$ sends $A(m)$ to $A(m-1)$ for $m\in [N_i+2,N_{i+1}+1]\cup [-N_{i+1}+1,-N_i]$ and $A(N_i+1)$ to $A(-N_i)$ and $A(-N_{i+1})$ to $A(N_{i+1}+1)$. Note also that $\Ad\,u_i$ fixes each element of  $A(m)$ for $m\in [-N_i+1,N_i]$. Namely $\Ad\,u_i\gamma$ acts trivially on
$$A[-N_{i+1},-N_i-1]\otimes A[N_i+1,N_{i+1}]
$$
and acts as a shift to the right on
$$A(-\infty,-N_{i+1}-1]\otimes A[-N_{i+1}+1,\infty)
$$
and cyclically on
$$
A[-N_i,N_i]
$$
respectively.
Let $\lambda=(\lambda_n)_{n\in\Z}$ be such that $\lambda_n=1$ for $n\in [-N_i, N_i]$ and $\lambda_n=0$ for $n\not\in [-N_{i+1},N_{i+1}]$ and let $H(\lambda)=\sum \lambda_n h_n$ where $h_n=h\in A(n)$. Since $\gamma(H(\lambda))=\sum \lambda_n h_{n+1}$, it follows that
$$
\Ad\,u_i\gamma(H(\lambda))=H(\lambda).
$$
Let $H_i=\sum_{-N_i\leq n\leq N_i}h_n$. We set $H_i(s)=(1-s)H_i+s H_{i+1}, s\in [0,1]$, which is a continuous path from $H_i$ to $H_{i+1}$ such that $\Ad\,u_i\gamma(H(s))=H(s)$.

Since $\Ad\,u_i\gamma(H_{i+1})=\Ad\,u_{i+1}\gamma(H_{i+1})$, it follows that $u_i^*u_{i+1}$ commutes with $\gamma(H_{i+1})$. By computation $\Ad(u_i^*u_{i+1})$ cyclically permutes the factors
$$A(m), \ m\in [-N_{i+1}+1, -N_i]\cup [N_i+1,N_{i+1}+1]
$$
upward and the factors
$$A(m), \ m\in[ -N_{i+2},-N_{i+1}]\cup [N_{i+1}+2,N_{i+2}+1]
$$
backward, respectively. Thus it follows that $u_i^*u_{i+1}$ is the tensor product of two matrices, each of which is a direct sum of permutation matrices. Regarding $u_i^*u_{i+1}$ as an element of
$$
B_i=\bigotimes_{m\in [-N_{i+2}, -N_i]\cup [N_i+1,N_{i+2}+1]}A(m)
$$
we obtain $\det(u_i^*u_{i+1})=1$ and a self-adjoint operator $L_i\in B_i$ such that $\|L_i\|\leq 1$, $\Tr(L_i)=0$, and $u_i^*u_{i+1}=e^{\pi i L_i}$. We define a path $s\mapsto u_i(s)=u_ie^{i\pi s L_i}$ from $u_i$ to $u_{i+1}$; its length is at most $\pi$ and $\det(u_i(s))=\det (u_i)$ and $\Ad\,u_i(s)\gamma(H_{i+1})=H_{i+1}$. One can apply the Rohlin property of $\gamma$ to the central sequence $(u_i(s))$ of unitary paths for a suitably chosen $(N_i)$. In the proof we invoke the following lemmas. (What we need for applying the Rohlin property is that any $\gamma$-cocycle of $s\mapsto u_i(s)$ in $C[0,1]\otimes A$ be connected to 1 by a unitary path of length at most a certain value, say $2\pi +1$.)

\begin{lem}
In the unitary group $\U=\U(C[0,1]\otimes M_n)$ let $\U_0$ be the set of $u\in \U$ such that each $u(s)$ has $n$ distinct eigenvalues for all $s\in[0,1]$. Then $\U_0$ is dense in $\U$. Moreover $\{u\in \U_0\ |\ \forall s\ \det(u(s))=1\}$ is dense in $\{u\in\U\ |\ \forall s\ \det(u(s))=1\}$.
\end{lem}
\begin{pf}
Probably this is well-known; see Lemma 3.3 of \cite{BEEK} for example. For any $u\in\U$ and $\epsilon>0$ there is a continuous function $h:[0,1]\ra (M_n)_{sa}$ such that $\|u(s)-e^{ih(s)}\|<\epsilon$. We may suppose that $h(s)$ has $n$ distinct eigenvalues. If $\det(u(s))=1$ then $e^{i\Tr(h(s))}\approx 1$ if $\epsilon$ is sufficiently small (depending on $n$). If $f(s)=\Tr(h(s))-2\pi k$ is small for some $k\in \Z$ we replace $h(s)$ by $h(s)-f(s)/n$. Then $s\mapsto v(s)=e^{ih(s)}$ approximates $u$ and satisfies that $\det(v(s))=1$.
\end{pf}

\begin{lem}
Let $u\in \U_0$ be such that $\det(u(s))=1$. Then there is a continuous function from $[0,1]$ into the self-adjoint elements of $M_n$ such that $u(s)=e^{ih(s)}$ and $\|h(s)\|<2\pi$ for $s\in [0,1]$. Hence $u$ can be connected to $1$ in $\U(A\otimes C[0,1])$ of length at most $2\pi$.
\end{lem}
\begin{pf}
Setting $h(0)=-i\log u(0)$ properly insures that $\Tr(h(0))=0$. Let $\lambda_{max}=\max \Spec(h(0))$ and $\lambda_{min}=\min \Spec(h(0))$. Since the gap on $\Spec(u(0))$ between $e^{i\lambda_{max}}$ and $e^{i\lambda_{min}}$  persists when $s$ moves from 0 to 1, one can continuously define $h(s)=-i\log u(s)$. Since $\Tr(h(s))\in \Z$ and $s\mapsto \Tr(h(s))$ is continuous it follows that $\Tr(h(s))=0$, which implies the gap never includes $1$, showing $\|h(s)\|<2\pi$. The last statement follows by defining a path $u_t,\ t\in [0,1]$ in $\U(A\otimes C[0,1])$ by $u_t(s)=e^{ith(s)}$.
\end{pf}

Let $k\in \N$ and let $u_i^{(k)}(s)=u_i(s)\gamma(u_i(s))\cdots \gamma^k(u_i(s))$, which is a $\gamma$-cocycle in $A\otimes C[0,1]$. What we need in applying the Rohlin property of $\gamma$ to a sequence $(u_i(s))$ is continuous paths in $\U(A\times C[0,1])$ connecting $u_i^{(k)}(s)$ with $1$, which form a central sequence as $i\ra\infty$ and $k=k_i\ra\infty$ for some $(k_i)$. Since $\det(u_i^{(k)}(s))$ is independent of $s$ and $u_i^{(k)}(s)$ is in a finite type I subfactor of $A$, which form a central sequence for a fixed $k$, we can find the desired paths for a sufficiently rapidly increasing $(N_i)$.

Thus we choose sequences $(z_i(s))$ and $(v_i(s))$ of paths in the unitary group of $A$ such that $u_i(s)=v_i(s)^*z_i(s)\gamma(v_i(s))$, $\max_s\|z_i(s)-1\|\ra0$, and $(v_i(s))$ is a central sequence. Let $H_i'(s)=v_i(s)H_i(s)v_i(s)^*$. Then $\Ad\,z_i(s)\gamma(H_i'(s))=H_i'(s)$. Note that
$$
u_{i+1}=v_i(1)^*z_i(1)\gamma(v_i(1))=v_{i+1}(0)^*z_{i+1}(0)\gamma(v_{i+1}(0)),
$$
and
$$
H_i'(1)=v_i(1) H_{i+1}v_i(1)^*\ \ {\rm and}\ \ H_{i+1}'(0)=v_{i+1}(0)H_{i+1}v_{i+1}(0)^*.
$$

\begin{lem}
For any $\epsilon>0$ and any finite subset $\F$ of $A$ there exists a $\delta>0$ and a finite subset $\G$ of $A$ with the following properties: Let $u,v_1,v_2$ be unitaries of $A$ such that $u,v_i,\gamma(v_i)$ are in a finite type I subfactor  of $A$ such that $\|v_i u\gamma(v_i^*)-1\|<\delta$ and $\|[x,v_i]\|<\delta, \ x\in\G$ for  $i=1,2$. Then for any $H=H^*\in A$ there are paths $\zeta(s), \ s\in[0,1]$ and $v(s),\ s\in [0,1]$ of unitaries in $A$ such that $\zeta(0)=1$, $\|\zeta(s) H\zeta(s)^*-H\|<\epsilon$, $v(0)=v_1$, $v(1)=v_2\zeta(1)$, $\|v(s)u\gamma(v(s)^*)-1\|<\epsilon$, and $\|[x,v(s)]\|<\epsilon$ and $\|[x,\zeta(s)]\|<\epsilon$ for any $x\in\F$.
\end{lem}
\begin{pf}
We shall indicate how to prove this lemma. Suppose that we have chosen $\delta>0$ and $\G$ with $\G^*=\G$.
Since $\|v_i^*\gamma(v_i)-u\|<2\delta$ it follows that $\|v_2v_1^*-\gamma(v_2v_1^*)\|<2\delta$. Since $\|[x,v_i]\|<\delta$ it also follows that $\|[x,v_2v_1^*]\|<2\delta,\ x\in \G$. That is, $v_2v_1^*$ is almost $\gamma$-invariant and central.

We define $k\in A_{sa}$ by functional calculus such that $\|k\|\approx0$ and
$$v_2v_1^*\gamma(v_2v_1^*)^*=e^{ik}.
$$
Let $\tau$ denote the unique tracial state of $A$. If $\tau(k)=0$ we set $\zeta(s)=1$; otherwise we choose a sufficiently central path $\zeta(s)$ of unitaries in a finite type I subfactor of $A$ (which commutes with the finite type I subfactor  containing $u,v_i,\gamma(v_i)$) such that $\gamma(\zeta(1))\approx \zeta(1)$ and  the self-adjoint $k'$ of small norm with $e^{ik'}=\zeta(1)\gamma(\zeta(1)^*)$ satisfies $\tau(k')=-\tau(k)$, which implies that the $k$ obtained by replacing $v_2$ by $v_2\zeta(1)$ above satisfies that $\tau(k)=0$. (The existence of such $\zeta(s)$ follows from the Rohlin property of $\gamma$; the centrality of $\zeta(s)$ is assured independently of $v_i$.)	Let $h\in A_{sa}$ be such that $v_2\zeta(1)v_1^*=e^{ih}$ and $\|h\|\leq 2\pi$ and remember that $e^{ik}=e^{ih} e^{-i\gamma(h)}$.

Let $w(s)=e^{ish}e^{-is \gamma(h)}e^{-isk}, \ s\in [0,1]$. Since $w(0)=w(1)=1$ we regard $w(s),s\in [0,1]$ as a unitary in $C(\T)\otimes A$ whose $K_1$ class is trivial as $\tau(k)=0$. Thus one can apply the Rohlin property to $w(s)$ to obtain a unitary $r\in C[0,1]\otimes A$ such that $r(0)=r(1)=1$, $w(s)\approx r(s)^*\gamma(r(s))$, or $\|r(s)e^{ish}-\gamma(r(s)e^{ish})e^{isk}\|\approx0$ and the centrality of $r(s)$ is determined by the one of $w(s)$ and the degree of the above approximation. Hence $w_1: s\in [0,1]\mapsto r(s)e^{ish}$ is an almost $\gamma$-invariant path from $1$ to $v_2\zeta(1)v_1^*$. Note that $v: s\mapsto w_1(s)v_1$ is a path from $w_1(0)v_1=v_1$ to $r(1)e^{ih}=v_2\zeta(1)v_1^*$ and $v(s)^*\gamma(v(s))=v_1^*w_1(s)^*\gamma(w_1(s)v_1)\approx v_1^*\gamma(v_1)\approx u$.
\end{pf}

By applying the above lemma we connect $H_i'(1)$ and $H_{i+1}'(0)$ as follows. For $u=u_{i+1}, v_1=v_i(1)$, and $v_2=v_{i+1}(0)$ we find $\zeta(s),v(s)$ as above; in particular $v(0)=v_i(1)$, $v(1)=v_{i+1}(0)\zeta(1)$. Then the path $v(s)H_{i+1}v(s)^*$ goes from
$$
H_i'(1)\ {\rm to}\ v_{i+1}(0)\zeta(1)H_{i+1}\zeta(1)^*v_{i+1}(0)^*
$$
which is connected by the inverse path of $v_{i+1}(0)  \zeta(s)  H_{i+1} {\zeta(s)}^*v_{i+1}(0)^*$ to
$$
 v_{i+1}(0)H_{i+1}v_{i+1}(0)^*=H_{i+1}'(0).
$$

Defining $z(s)\approx1$ by $v(s)^*z(s)\gamma(v(s))=u_{i+1}$ we have that $\Ad\, z(s)\gamma(v(s)H_{i+1}v(s)^*)=v(s)H_{i+1}v(s)^*$, where $z(0)=v_1u_{i+1}\gamma(v_1)^*=z_i(1)$ and $z(1)=v_2\zeta(1)u_{i+1}\gamma(v_2\zeta(1))^*$. Since we can assume that $\zeta(s)$ can be arbitrarily central we may assume that $\zeta(s)H_{i+1}\zeta(s)^*\approx H_{i+1}$ with an arbitrary precision. Hence we may assume  that
$$
\Ad\, z(1)\gamma(v_2\zeta(s)H_{i+1}\zeta(s)^*v_2^*)\approx v_2\zeta(s)H_{i+1}\zeta(s)^*v_2^*.
$$
Then we have, at $s=0$, the end of path, that
$$
\Ad\, z(1)\gamma(H_{i+1}'(0))\approx H_{i+1}'(0)=\Ad\, z_{i+1}(0)\gamma(H_{i+1}'(0)).
$$
Since both $z(1)$ and $z_{i+1}(0)$ are close to $1$ we find $h=h^*$ such that $z_{i+1}(0)z(1)^*=e^{ith}$ and $\|[h,\gamma(H_{i+1}'(0))]\|\approx0$. Thus the path $w: s\mapsto e^{ish}z(1)$ in a small vicinity of 1 satisfies that $\Ad\,v(s)\gamma(H_{i+1}'(0))\approx H_{i+1}'(0)$.  Thus by connecting the paths $H_i'(s)$ we have obtained a long path $H(s)$ such that $\alpha_t(x)=\lim_{s\ra\infty} \Ad\,e^{itH(s)}(x)$ and $\|\Ad\,z(s)\gamma(H(s))-H(s)\|\ra0$ for some function $z:s\mapsto \U(A)$ with $\|z(s)-1\|\ra0$. (We have chosen $z(s)$ to be continuous; but this is not necessary and can be assumed automatically as above.) This concludes the proof of Proposition \ref{no-interaction}.

\begin{theo}
Let $\beta$ be a flow on a UHF algebra $F$ and define a flow $\alpha$ on $A=\bigotimes_{m\in\Z}A(m)$ with $A(m)=F$ by $\alpha_t=\bigotimes_{m\in\Z} \beta_t$. Let $\bar{\alpha}$ denote the extension of $\alpha$ to a flow on the crossed product $A\times_\gamma\Z$. If $\alpha$ is continuously AI then so is $\bar{\alpha}$.
\end{theo}

Since $F$ is a UHF algebra there is an increasing sequence $F_n$ of C$^*$-subalgebras such that $1_A\in F_n$, $F_n\cong M_{k_n}$ for some $k_n\in \N$, and the union of $F_n$ is dense in $A$. Since $\beta$ is continuously AI there is a continuous function $h$ from $[0,\infty)$ into $A_{sa}$ such that $\beta_t=\lim_{s\ra\infty}\Ad\,e^{ith(s)}$ uniformly in $t\in [-1,1]$. By changing $h$ slightly and reparameterizing $s$ we may suppose that $h(s)\in F_n$ for $s\leq n$. We also assume that $h(0)=0$.

To prove that $\bar{\alpha}$ is continuously AI we recall how we proved it in the case $h(s)$ is just constant $h$. Let $H_i=\sum_{-N_i\leq m\leq N_i} h_m$ where $h_m$ is $h$ in $A(m)$ for $m\in \Z$ and $(N_i)$ is a certain increasing sequence in $\N$. Define $H(s)=H_i+(s-i)(H_{i+1}-H_i)=(i+1-s)H_i+(s-i)H_{i+1}$ when $i\leq s\leq i+1$. Then $\alpha_t=\lim_{s\ra\infty}\Ad\,e^{itH(s)}$ on $A=\bigotimes_{\Z}A(m)$. We have chosen a sequence $(u_i)$ and $(v_i), (z_i)$ such that $z_i\approx1$ and $u_i= v_i^*z_i\gamma(v_i)$ (associated with $N_i$) and defined $H'(s)=v_iH(s)v_i^*, \ i\leq s<i+1$. Then $H'(s)$ is not continuous at $s=i+1$. To remedy this we have chosen an appropriate path $u_i(s), 0\leq s\leq 1$ from $u_i$ to $u_{i+1}$ and obtain $v_i(s),z_i(s)\approx 1$ such that $u_i(s)=v_i(s)^*z_i(s)\gamma(v_i(s))$ with $v_i(0)=v_i$ and $z_i(0)=z_i$ and define a path $H_1(s)=v_i(s)H_{i+1}v_i(s)^*$ from $H'(i+1-0)$. But $H_1(1)$ is still not equal to $H'(i+1+0)$ because we did not place the condition $v_i(1)=v_{i+1}$ (though $v_i(1)^*\gamma(v_i(1))\approx u_{i+1}\approx v_{i+1}^*\gamma(v_{i+1})$. (What we actually did was a bit different from this; we combined these two processes.) So we have chosen anther path $v(s)$ from $v_i(1)$ to $v_{i+1}\zeta(1)$ where $\zeta(s)$ is a sufficiently central path starting from 1. We have defined $H_2(s)=v(s)H_{i+1}v(s)^*$ and $H_3(s)=v_{i+1}\zeta(1-s)H_{i+1}\zeta(1-s)^*v_{i+1}^*$. We have inserted $H_1(s),H_2(s),H_3(s)$ between $H'(i+1-0)$ and $H'(i+1+0)$.

In the present case $u_i$ will depend not only on $N_i,N_{i+1}$ but also on $F_n$, where in particular $\Ad\,u_i^{(n)}\gamma$ should be trivial on
$$
F_n[-N_{i+1},-N_i-1]\otimes F_n[N_i+1,N_{i+1}];
$$
so let us write $u_i^{(n)}$ instead of $u_i$. We define
$$
H(s)=\sum_{|m|\leq N_i}h_m(s)+\sum_{N_i<|m|\leq N_{i+1}}h_m((i+1)(s-i)) \in \bigotimes_{-N_{i+1}}^{N_{i+1}}F_i
$$
for $i\leq s\leq i+1$; then it follows that $\Ad\,u_i^{(i)}\gamma(H(s))=H(s)$ for $i\leq s\leq i+1$. Then by taking $u_i^{(i)}$ instead of $u_i$ we can proceed just as before. What we really need is $\det((u_i^{(i)})^*u_{i+1}^{(i+1)})=1$ in $F_{i+1}[-N_{i+2},-N_i]\otimes F_{i+1}[N_i+1,N_{i+2}+1]$, which follows as before. We will then obtain a self-adjoint $L_i\in F_{i+1}[-N_{i+2},-N_i]\otimes F_{i+1}[N_i+1,N_{i+2}+1]$ such that $(u_i^{(i)})^*u_{i+1}^{(i+1)}=e^{i\pi L_i}$ and $\|L_i\|\leq 2$.

%%%%%%%%%%%%%%%%%%%%%%%%%%%%%%%%%%%%%%%%%%%%%%%%%%%%%%%%%%%%%%%%%%%%%%%%

\section{Concluding Remarks}

Let $B$ be a unital nuclear C$^*$-algebra and let $A(k)=B$ for all $k\in \Z$. If $K_0(B)$ is not of rank 1 then the shift automorphism $\gamma$ on $A=\bigotimes_{k\in \Z} A(k)$ is not approximately inner (because $\gamma$ acts on $K_0(A)$ in a non-trivial way). If $[1]$ is not divisible by any positive integer $>1$ in $K_0(B)$ then $\gamma$ does not have the {\em cyclic} Rohlin property which was used in a crucial way for proving $\bar{\alpha}$ is AI. Let $\Phi$ be a $\gamma$-invariant potential in $A$ such that $\|\Phi\|_\lambda<\infty$ for some $\lambda>0$ or $\Phi(\Lambda)$'s commute with each other. Then, as in the case $B$ is finite-dimensional, one can define a flow $\alpha^\Phi$ on $A$ based on $\Phi$.

\begin{prop}
In the above situation suppose that $B$ is quasi-diagonal. Then the flow $\bar{\alpha}^\Phi$ is quasi-diagonal on $A\times_\gamma\Z$.
\end{prop}
\begin{pf}
We may assume that $\Phi$ has finite range.
For a positive integer $n$ let $A_n=\bigotimes_{k\in \Z_{2n+1}}A(k)$, the tensor product of $2n+1$ copies of $B$. Let $\Phi_n$ be the potential obtained from $\Phi$ by imposing the periodic boundary condition on $\{-n,-n+1,\ldots,0,1,2,\ldots,n-1\}$. The flow $\alpha^{\Phi_n}$ defined on $A_n$ is inner and so is MF (as $A_n$ is MF). The same is true for the flow $\bar{\alpha}^{\Phi_n}$ on $A_n\times_{\gamma_n}\Z_{2n+1}$ (as $A_n\times_{\gamma_n}\Z_{2n+1}$ is MF), where $\gamma_n$ is the natural action of $\Z_{2n+1}$ on $A_n$. We obtain a continuous field of flows on $\N \cup \{\infty\}$ by associating $\bar{\alpha}^{\Phi_n}$ to $n\in \N$ and $\bar{\alpha}$ to $\infty$ by associating continuous fields of operators with $C_{00}(\Z,A)$. (Here $A_n\times_{\gamma_n}\Z_{2n+1}$ is identified with $C(\{-n,-n+1,\ldots,n\},\bigotimes_{|k|\leq n}A(k))\subset C_0(\Z,A)$. See the proof of Proposition 4.1 of \cite{AI-MF} for details.) Since the flows over $\N$ are MF it follows from Proposition 3.8 of \cite{AI-MF} that the flow on $\infty$ is MF. Since $A\times_\gamma\Z$ is nuclear we conclude that $\bar{\alpha}^\Phi$ is QD.
\end{pf}

Our aim was to show two things; The flows could be QD (quasi-diagonal) without being AI (approximately inner) and AI without being continuously AI (or asymptotically approximately inner) for simple (or technically simple) C$^*$-algebras. To achieve them we considered an extension of an AI flow to a crossed product. But we failed because we could not come up with a new way of disproving AI or continuous AI and ended up to give a new class of AI flows. There is still a possibility that the approach may not be futile.

There was one more thing we wanted to prove; non-AI flows (or non-MF flows) could satisfy the KMS condition (the existence of KMS states for all inverse temperatures) for a simple C$^*$-algebra, which still eluded our grasp. We presented a non-simple (even non-prime) example based on a flow on the Toeplitz algebra.

\end{document}